\documentclass[review,hidelinks,onefignum,onetabnum]{siamart220329}

\usepackage{lipsum}
\usepackage{amsfonts}
\usepackage{graphicx}
\usepackage{epstopdf}
\usepackage{algorithmic}
\usepackage{esint}
\usepackage{subfigure}
\usepackage{multirow}
\usepackage[export]{adjustbox}
\ifpdf
  \DeclareGraphicsExtensions{.eps,.pdf,.png,.jpg}
\else
  \DeclareGraphicsExtensions{.eps}
\fi

\newsiamremark{remark}{Remark}
\newsiamremark{hypothesis}{Hypothesis}
\crefname{hypothesis}{Hypothesis}{Hypotheses}
\newsiamthm{claim}{Claim}

\usepackage{amsfonts,amssymb,amsmath}
\newcommand{\tcb}[1]{\textcolor{blue}{#1}}

\newcommand{\E}{\mathbb{E}}

\def\b{\boldsymbol}

\def\e{\epsilon}

\usepackage{amsopn}

\numberwithin{equation}{section}
\newtheorem{example}{Example}[section]

\title{Random ordinate method for mitigating the ray effect in radiative transport equation simulations\thanks{\today
\funding{Lei Li is partially supported by the National Key R\&D Program of China No. 2020YFA0712000; Shanghai Municipal Science and Technology Major Project 2021SHZDZX0102, NSFC 12371400 and NSFC 1203101, the Strategic Priority Research Program of Chinese Academy of Sciences, Grant No. XDA25010403. Min Tang is funded by The Strategic Priority Research Program of Chinese Academy of Sciences, No.XDA25010401, NSFC12411530067, NSFC12031013, and partially by Mevion Medical Equipment Co., Ltd.}}}

\author{Lei Li\thanks{School of Mathematical Sciences, Institute of Natural Sciences, MOE-LSC,
Shanghai Jiao Tong University, Shanghai, P.R. China.
  (\email{leili2010@sjtu.edu.cn} 
  ).}
\and Min Tang\thanks{School of Mathematical Sciences,, Institute of Natural Sciences and MOE-LSC,
Shanghai Jiao Tong University, Shanghai, P.R. China. 
  (\email{tangmin@sjtu.edu.cn}).}
\and Yuqi Yang\thanks{School of Mathematical Sciences,, Institute of Natural Sciences, Shanghai Jiao Tong University, Shanghai, P.R. China. 
  (\email{yangyq1@sjtu.edu.cn}).}}
\ifpdf
\hypersetup{
  pdftitle={},
  pdfauthor={}
}
\fi

\begin{document}

\maketitle

\begin{abstract}
The Discrete Ordinates Method (DOM) is the most widely used velocity discretization method for simulating the radiative transport equation. However, the ray effect is a long-standing drawback of DOM. In benchmark tests that exhibit the ray effect, we observe low regularity in the velocity variable of the solution. To address this issue, we propose a Random Ordinate Method (ROM) to mitigate the ray effect. Compared to other strategies proposed in the literature for mitigating the ray effect, ROM offers several advantages: 1) For benchmark tests that exhibit ray effect, the computational cost is lower than that of the DOM; 2) it is simple and requires minimal changes to existing DOM-based code; 3) it is easily parallelizable and independent of the problem setup. A formal analysis is presented for the convergence orders of the error and bias. Numerical tests demonstrate the reduction in computational cost compared to DOM, as well as its effectiveness in mitigating the ray effect.
\end{abstract}



\section{Introduction}
The radiative transport equation (RTE) stands as a fundamental equation governing the evolution of angular flux as particles traverse through a material medium.  
It provides a statistical description of the density distribution of particles. The RTE has found extensive applications across diverse fields, including astrophysics \cite{1999Hydrodynamics, 2008Manufactured}, fusion \cite{2005Pulsed, 2001Three}, biomedical optics \cite{2010PDE, 2008Light}, and biology, among others.

The steady state RTE with anisotropic scattering reads
\begin{equation}\label{eq:equ_1}
\b{u} \cdot \nabla \psi(\b{z}, \b{u})+\sigma_{T}(\b{z}) \psi(\b{z}, \b{u})=\sigma_{S}(\b{z})\frac{1}{|S|}\int_{S} P(\b{u'},\b{u})\psi(\b{z}, \b{u'}) \mathrm{d} \b{u'}+q(\b{z}),
\end{equation}
subject to the following inflow boundary conditions:
\begin{equation}\label{eq:equ_1_bc}
	\psi(\boldsymbol{z}, \boldsymbol{u})=\psi_{\Gamma}^{-}(\boldsymbol{z}, \boldsymbol{u}), \quad \boldsymbol{z} \in \Gamma^{-}=\partial \Omega, \quad \boldsymbol{u} \cdot \boldsymbol{n}_{\boldsymbol{z}}<0.
\end{equation}
Here, $\boldsymbol{z} \in \Omega \subset \mathbb{R}^3$ represents the spatial variable; $\boldsymbol{u}$ denotes the direction of particle movement, and $S=\{\boldsymbol{u} \mid \boldsymbol{u} \in \mathbb{R}^3, |\boldsymbol{u}|=1 \}$; $\boldsymbol{n}_{\boldsymbol{z}}$ stands for the outward normal vector at position $\boldsymbol{z}$. $\psi(\boldsymbol{z}, \boldsymbol{u})$ gives the density of particles moving in the direction $\boldsymbol{u}$ at position $\boldsymbol{z}$. The coefficients $\sigma_{T}(\boldsymbol{z})$, $\sigma_{S}(\boldsymbol{z})$, and $q(\boldsymbol{z})$ represent the total, scattering cross-sections, and the source term, respectively. The cross-sections are bounded and $ 0\le \sigma_{S}(\boldsymbol{z}) <\sigma_{T}(\boldsymbol{z})<\infty$, the source term $q(\boldsymbol{z}) \in L^{\infty}(\Omega)$ is non-negative. For physically meaningful situations, $\sigma_{S}(\boldsymbol{z})/\sigma_{T}(\boldsymbol{z})\leq 1$, for $\forall \boldsymbol{z}\in \Omega$,  while in this paper we consider the case when $\sigma_{S}(\boldsymbol{z})/\sigma_{T}(\boldsymbol{z})\leq \lambda<1$. The kernel $k(\b{u}', \b{u}):=\frac{1}{|S|}P(\b{u}',\b{u})$ is the scattering kernel, which provides the probability that particles moving in the direction $\b{u'}$ scatter to the direction $\b{u}$.
For the notational convenience, we will use the symbol $\fint_S:=\frac{1}{|S|}\int_S$ to denote the average over the domain $S$ associated with the indicated measure. Then, the scaled scattering kernel satisfies:
$$P(\b{u}',\b{u})=P(\b{u},\b{u'}),\qquad \fint_S P(\b{u'},\b{u})\mathrm{d} \b{u} =\frac{1}{|S|}\int_SP(\b{u'},\b{u})\mathrm{d} \b{u}=1.$$

The numerical methods for solving the RTE are mainly divided into two categories: particle methods and PDE-based methods. Particle methods, specifically Monte Carlo (MC) methods, simulate the trajectories of numerous particles and collect the density distribution of all particles in the phase space to obtain the RTE solution. The MC methods are known to be slow and noisy but are easy to parallelize and suitable for all geometries \cite{1971An, bhan2007condensed}. Meanwhile, the PDE-based methods are more accurate and can be faster, but they are not as flexible as the MC method for parallel computation and complex geometries \cite{lewis1984computational}. In this paper, we are interested in the PDE-based method.

The Discrete Ordinates Method (DOM) \cite{1970Nuclear, Carlson1955SOLUTION} is the most popular velocity discretization method. DOM approximates the solution to \eqref{eq:equ_1_bc} using a set of discrete velocity directions $\b{u_m}$, which are referred to as ordinates. The integral term on the right-hand side of Equation \eqref{eq:equ_1} is represented by weighted summations of the discrete velocities. DOM retains the positive angular flux and facilitates the determination of boundary conditions. However, solving the RTE using the standard DOM is expensive due to its high dimensionality since the RTE has three spatial variables and two velocity direction variables.

In real applications, people are usually interested in the spatial distributions of some macroscopic quantities, such as the particle density $\int_S\psi(\b{z}, \b{u}) \mathrm{d}  \b{u}$, the momentum $\int_S \b{u}\psi( \b{z},  \b{u})\mathrm{d} \b{u}$, etc. Thus, the spatial resolution must be adequate. The computational costs can be significantly reduced if one can obtain the right macroscopic quantities by using a small number of ordinates in DOM.

One natural question arises: Can we improve the accuracy of macroscopic quantities without increasing the number of ordinates? Considerable effort has been invested in discovering a quadrature set with high-order convergence. For the 1D velocity in slab geometry, Gaussian quadrature exists, exhibiting spectral convergence when the solution maintains sufficient smoothness in velocity. However, devising a spectrally convergent 2 or 3-dimensional Gauss quadrature remains unclear\cite{Larsen2010}.Furthermore, in section 2, we highlight, through numerical tests, that the solution's regularity in the velocity direction can be considerably low in some benchmark tests. Consequently, it remains uncertain whether better approximations can be expected, even if a 2 or 3-dimensional Gaussian quadrature with high convergence order for smooth solutions is identified.
 


When employing a limited number of ordinates in DOM, the ray effect becomes noticeable in numerous 2D spatial benchmark tests \cite{2013MomentClosures, LathropRayeffect}. The macroscopic particle density exhibits nonphysical oscillations along the ray paths, particularly noticeable when the inflow boundary conditions or radiation sources demonstrate strong spatial variations or discontinuities. The ray effect stems from particles being confined to move in a limited number of directions. As highlighted in section 2, in benchmark tests displaying the ray effect, we observe low regularity in velocity within the solution, indicating a low convergence order for DOM. In order to mitigate the ray effect, one has to increase the number of ordinates \cite{1977Ray}, which will significantly increase the computational cost.

The ray effect stands as a long-standing drawback in DOM simulations, and several strategies have been proposed to mitigate these ray effects at reasonable costs. Examples include approaches discussed in \cite{Ray_effect_mitigation_for_the_discrete_ordinates_method_through_quadrature_rotation, Artificial_Scattering, MorelAnalysisofRayEffect, Frame_Rotation}. The main idea is to use biased or rotated quadratures or combine the spectral method with DOM. In this paper, inspired by the randomized integration method \cite{Novak1988}, we introduce a random ordinate method (ROM) for solving RTE. Compared with other ray effect mitigating strategies proposed in the literature, the advantages of ROM are: 1) the computational cost is comparable to DOM; 2) it is simple and makes almost no change to all previous code based on DOM; 3) It is easy to parallelize and independent of the problem setup.

Randomized algorithms can achieve higher convergence order when the solution regularity is low \cite{Novak1988}. The concept of the randomized method is straightforward: when approximating $\int_{0}^1 f(x) dx$, the integral interval $[0,1]$ is partitioned into cells with a maximum size of $h$. Subsequently, an $x_m$ is chosen from each interval, and $\int_{0}^1 f(x) dx$ is approximated by $\sum_{\ell=1}^n \omega_{\ell} f(x_{\ell})$, where $\omega_{\ell}$ represents the quadrature weights. With a fixed set of $\{x_{\ell},\omega_{\ell}\}$, one can expect uniform first-order convergence for general $f(x)$ that is Lipschitz continuous. On the other hand, when one randomly chooses a point $x_{\ell}$ inside each interval and keeps using the same $\omega_{\ell}$, the expected error can achieve $O(h^\frac{3}{2})$ convergence, whereas the expectation of all randomly chosen quadrature provides $O(h^3)$ convergence. Randomized algorithms have been developed and analyzed for various problems, such as the initial value problem of ordinary differential equation (ODE) systems \cite{ErroranalysisSTENGL1995,STENGLE199025} whose complexity has been studied in \cite{HEINRICH200877, 2011On},  for stochastic differential equations \cite{Clark1980} whose complexity is analyzed in \cite{Cao2021}, and the interacting particle system \cite{jin2020random,jin2021convergence}.

The main idea of ROM is that the velocity space is partitioned into $n$ cells, and a random ordinate is selected from each cell. A DOM system with those randomly chosen ordinates is then solved. ROM runs a lot of samples with $n$ ordinates and then we calculatetheir expectation. When the solution regularity is low in the velocity variable, even with a carefully chosen quadrature in DOM, it can only achieve a low convergence order, which is similar as randomly chosen ordinates. According to the formal analysis and numerical results, the ROM solutions' expected values can achieve a higher convergence order in the velocity space. Consequently, the accuracy doesn't decrease for a single run, while averaging multiple runs can lead to a higher-order convergence and mitigate the ray effect. Since different runs employ distinct ordinates, ROM allows for easy parallel computation.

The novelty of the ROM lies in its capability to employ deterministic solvers with lower convergence orders and obtain the numerical results with higher order convergence by leveraging the expectation of random sampling. The ROM advances the state of the art by offering a method that is both computationally efficient—with a cost comparable to DOM—and exceptionally straightforward to implement, requiring minimal modifications to existing DOM codes. Crucially, its independence from problem geometry and inherent parallelizability make it highly attractive for complex practical applications. Extensive numerical results validate its effectiveness, demonstrating superior convergence order and a dramatic reduction in ray effects compared to the standard DOM at equivalent computational costs.



This paper is structured as follows: Section 2 delves into the ray effects of DOM and illustrates the low regularity of the solution in velocity space. Details of the ROM are described in Section 3. In Section 4, analytical results are presented for the convergence orders of the error and bias when ROM is applied to RTE with isotropic scattering in slab geometry. Section 5 displays the numerical performance of the ROM, demonstrating its ability to mitigate the ray effect. Finally, Section 6 concludes the paper with discussions.

\section{Ray effects and low regularity}\label{Motivation}

The DOM is the most popular angular discretization for RTE simulations, it writes \cite{Larsen2010}: 
$$
\b{u}_{\ell} \cdot \nabla \psi_{\ell}(\b{z} )+\sigma_{T}(\b{z}) \psi_{\ell}(\b{z})=\sigma_{S}(\b{z}) \sum_{\ell'\in V} w_{\ell'}P_{\ell,\ell'}\psi_{\ell'}(\b{z})+q_{\ell}(\b{z}),\quad \ell\in V,
$$
subject to the following inflow boundary conditions:
$$
	\psi(\b{z}, \b{u}_{\ell})=\psi_{\ell}^{-}(\b{z}, \b{u}), \quad \b{z} \in \Gamma^{-}=\partial \Omega, \quad \b{u}_{\ell} \cdot \b{n}_{\b{z}}<0.
$$
where $\omega_{\ell}$ is the weight of the quadrature node $\b{u}_{\ell}$, satisfying  $\sum_{\ell \in V} \omega_{\ell}=1$; $V$ represents the index set of the quadrature $\{\b{u}_{\ell},\omega_{\ell}\}$;  $P_{\ell,\ell'}\approx P(\b{u_{\ell}},\b{u_{\ell'}})$ denotes the discrete scattering kernel; $q_{\ell}(\b{z})\approx q\left(\b{z}, \b{u}_{\ell}\right)$. Then $\psi_{\ell}(\b{z}) \approx \psi\left(\b{z},\b{u}_{\ell}\right)$ and, for all $\ell \in V$,
$$\sum_{\ell'\in V} \omega_{\ell'}P_{\ell,\ell'}\psi_{\ell'}(\b{z})\approx \fint_{S} P(\b{u_{\ell}},\b{u}_{\ell'})\psi(\b{z}, \b{u}_{\ell'}) \mathrm{d} \b{u}_{\ell'}.$$

Two commonly used quadrature methods are Uniform Quadrature and Gaussian Quadrature. For detailed information, please refer to Appendix \ref{app:detailsDOM}.

It has long been known that the solution of DOM exhibits the ray effect, especially when there are discontinuous source terms in the computational domain \cite{Ray_effect_mitigation_for_the_discrete_ordinates_method_through_quadrature_rotation, Artificial_Scattering, Frame_Rotation}. This phenomenon cannot be improved by increasing the spatial resolution. 
We show one typical example  to demonstrate the ray effects.
\begin{example}\label{2D disk_ray effect}
We consider RTE in the X-Y geometry with a localized source term at the center of the computational domain. Let
$$
x\times y \in \Omega =[0,1] \times [0,1],\quad
\sigma_T=1,\quad \sigma_S=0.5.
$$
$$ 
q(x,y)=\left\{
	\begin{aligned}
	2, & \quad (x,y)\in[0.4,0.6]\times [0.4,0.6],\\
	0, & \quad \mbox{elsewhere}.\\
	\end{aligned}
	\right	.
$$
The inflow boundary conditions are zero.
\end{example}

We consider the isotropic and anisotropic scattering with the following scattering kernel 
\begin{equation}
P(\b{u},\b{u'})=G(\b{u} \cdot \b{u'})=G(\cos \xi)=1+g \cdot \cos\xi,
\label{eq:anisotropic kernel}
\end{equation}
where $\xi$ is the included angle between $\b{u}$ and $\b{u'}$. When $g=0$, \eqref{eq:anisotropic kernel} gives isotropic scattering, meaning particles moving with velocity $\b{u}$ will scatter into a new velocity $\b{u}'$ with uniform probability for all $\b{u}'$. When $g=0.9$, \eqref{eq:anisotropic kernel} gives anisotropic scattering, implying that when the included angle $\xi$ between $\b{u}$ and $\b{u'}$ is smaller, the probability that particles moving with velocity $\b{u}$ scatter into $\b{u}'$ is higher.

We would like to emphasize that the main purpose of our current work is to show the convergence orders of the velocity discretization; any spatial discretization can be chosen to obtain the numerical results. Since only the convergence order in the velocity variable is considered in this paper, after spatial discretization, the reference solution can be considered as a large integral system for unknowns at the given spatial grids. We use the same spatial mesh and discretizations throughout the paper, eliminating the need to consider the error introduced by the spatial discretization. 

In the paper, the classical second-order diamond difference (DD) method \cite{Edward1989Asymptotic,lewis1984computational} is adopted for the spatial discretization to solve Example \ref{2D disk_ray effect}. The spatial domain $\Omega=[x_L,x_R]\times[y_B,y_T]$ is divided into  $I \times J$  uniform cells. Let 
$\Delta x=\frac{x_R-x_L}{I}$, $\Delta y=\frac{y_T-y_B}{J}$, $x_0=x_L$, $y_0=y_B$ and
\begin{align*}
&x_i=x_L+i\Delta x,\quad x_{i+\frac{1}{2}}=x_L+\big(i-\frac{1}{2}\big)\Delta x,
\quad\mbox{for $i=1,\cdots,I$,}\\
&y_j=y_B+j\Delta y,\quad
 y_{j+\frac{1}{2}}=y_B+\big(j-\frac{1}{2}\big)\Delta y, \quad\mbox{for $j=1,\cdots,J$.}
\end{align*}

We use DD with $I=100$ and $J=100$. The grid points of the two-dimensional DD method are at the cell centers; that is, approximations of $\psi_\ell(x_{i-\frac{1}{2}},y_{j-\frac{1}{2}})$ and the average density 
\[
\phi\left(x_{i-\frac{1}{2}},y_{j-\frac{1}{2}}\right) = \sum_{\ell \in V}\bar{\omega}_{\ell}\psi_{\ell}\left(x_{i-\frac{1}{2}},y_{j-\frac{1}{2}}\right)
\]
(for $i=1,\cdots,I$, $j=1,\cdots,J$) are obtained, where $\bar{\omega}_{\ell}$ are weights in X-Y geometry.

\begin{figure}[htbp]
	\centering

\subfigure[4 ordinates ]
{
		\includegraphics[width=0.3\linewidth]{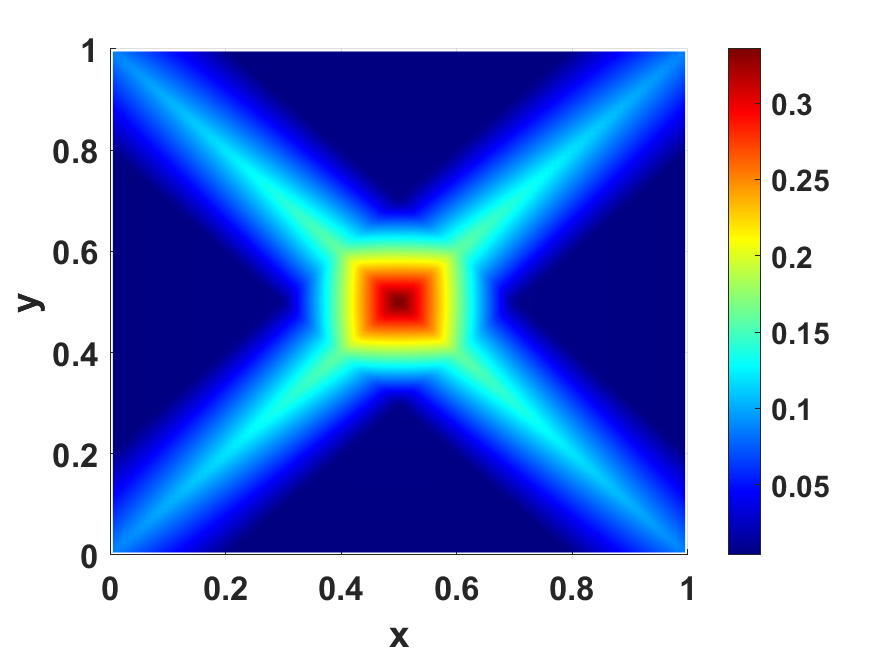}
}
\subfigure[16 ordinates ]
{
		\includegraphics[width=0.3\linewidth]{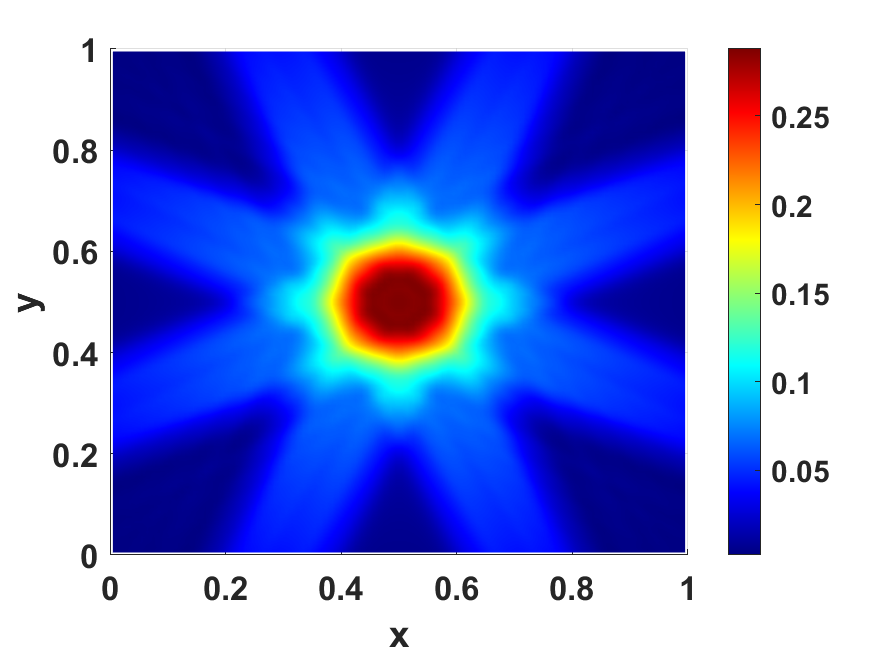}
}
\subfigure[10000 ordinates ]
{

		\includegraphics[width=0.3\linewidth]{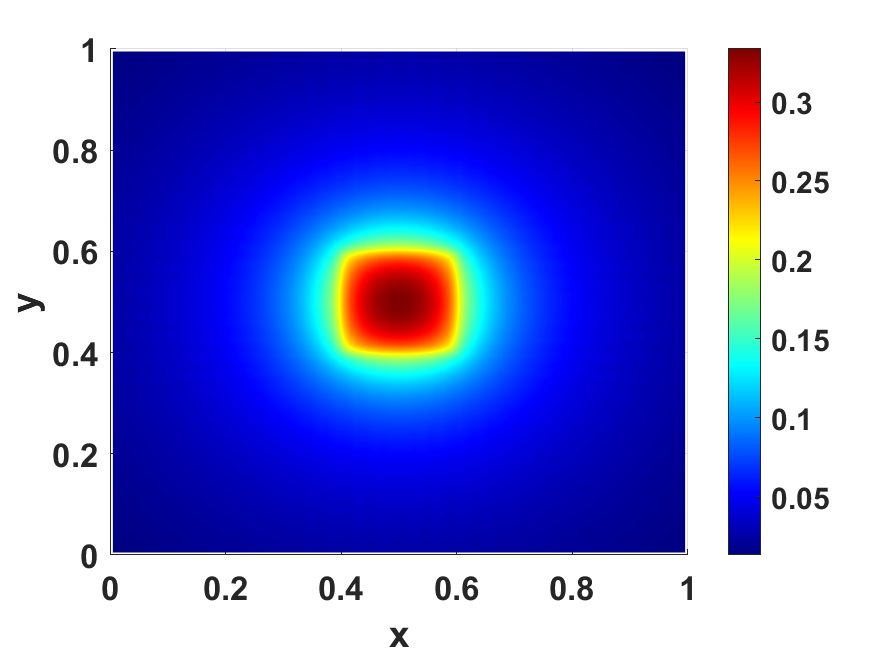}
}

\subfigure[4 ordinates ]
{
		\includegraphics[width=0.3\linewidth]{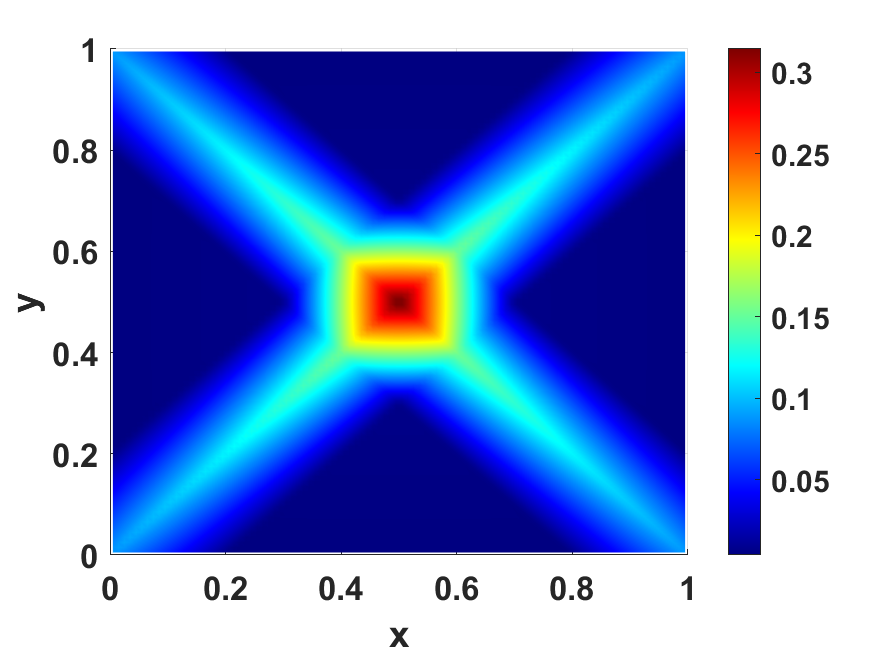}
}
\subfigure[16 ordinates ]
{
		\includegraphics[width=0.3\linewidth]{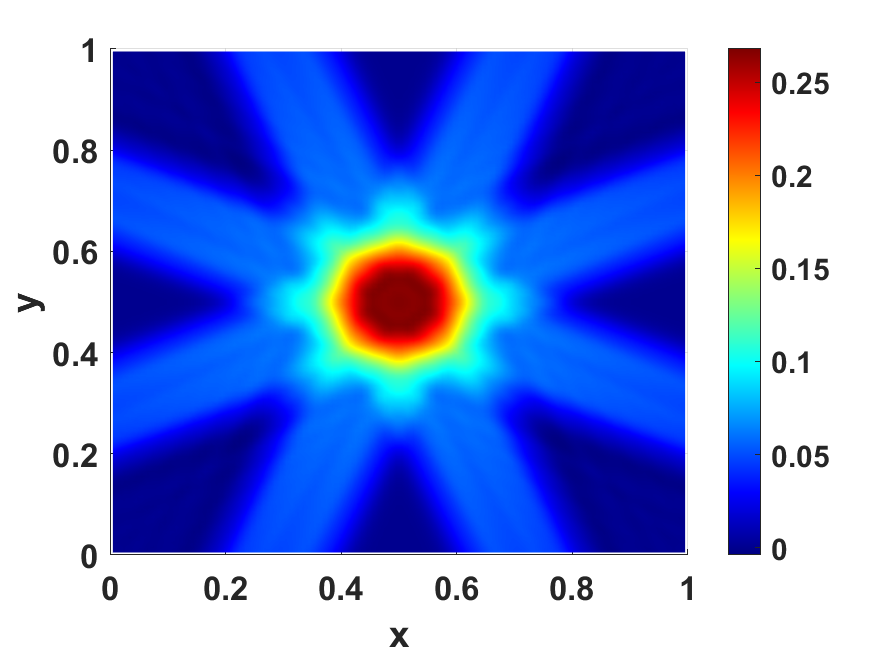}
}
\subfigure[10000 ordinates ]
{

		\includegraphics[width=0.3\linewidth]{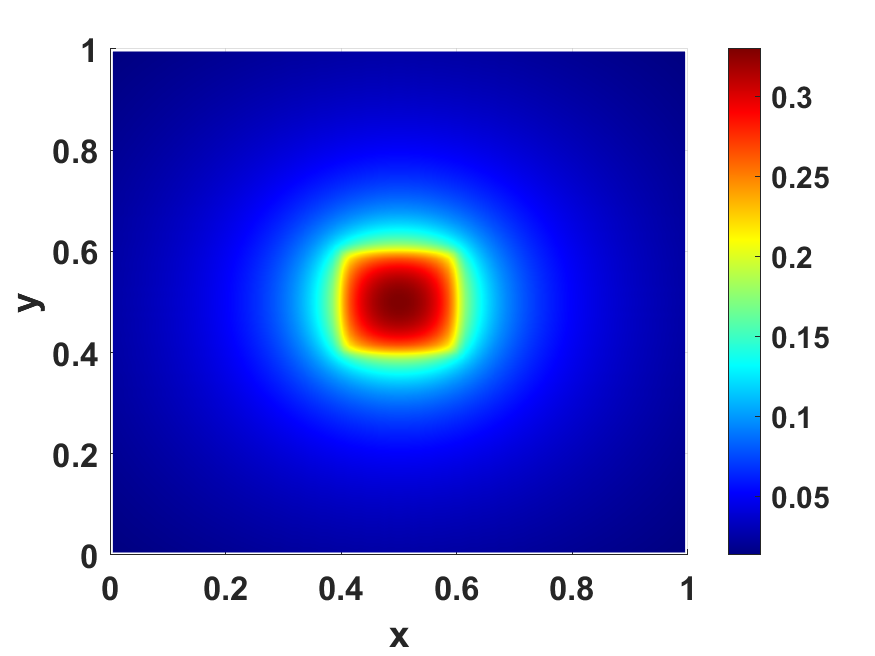}
}

		\caption{Demonstration of the ray effects. The average densities $\phi$ calculated with different numbers of ordinates are displayed. The numerical results are calculated with $100 \times 100$ spatial cells and uniform quadrature. (a)(b)(c): isotropic scattering kernel with $g=0$ in \eqref{eq:anisotropic kernel}; (d)(e)(f) anisotropic scattering kernel with $g=0.9$ in \eqref{eq:anisotropic kernel}.}
		\label{fig:section2_fig_rayeffect}
\end{figure}

The numerical results of the average density $\phi$ with various numbers of ordinates are displayed in Figure \ref{fig:section2_fig_rayeffect}. From left to right, $4$, $16$, and $10000$ ordinates of Uniform quadrature 
. We can observe that the solutions exhibit rays that correspond to the chosen ordinates of the DOM. The solutions have poor accuracy and this phenomenon does not disappear with the refinement of the spatial grids.

The \(\ell^2\) errors between the reference solution \(\phi^{ref}\) and the numerical solutions \(\phi\) obtained by different quadratures are defined by

\begin{equation}\label{eq:l2norm2D}
\mathcal{E}_{DOM}=\sqrt{\frac{1}{I J}\Big(\sum_{i=0}^{I-1}\sum_{j=0}^{J-1}\mid\phi(x_{i+\frac{1}{2}},y_{j+\frac{1}{2}})-\phi^{ref}(x_{i+\frac{1}{2}},y_{j+\frac{1}{2}})\mid^{2}\Big)}.
\end{equation}

The reference solution of the Uniform (Gaussian) quadrature \(\phi_{U}^{ref}\) (\(\phi_{G}^{ref}\)) is computed by \(N=20\), which indicates \(1600\) (\(840\)) ordinates on the 2D disk. As shown in Table \ref{tabel:2D_DOM_order}, the convergence orders for both quadrature methods are approximately $0.8$. This is a consequence of the solution's low regularity in velocity space. It is well-known that improving solution accuracy is challenging without increasing the number of quadrature nodes under conditions of low regularity.

\begin{table}[!ht]
    \centering
    \caption{Example \ref{2D disk_ray effect}: Convergence orders of DOM with the scattering kernel as in \eqref{eq:anisotropic kernel}. Here $\Delta S=\frac{\pi}{4M}$.}
    \label{tabel:2D_DOM_order}
    \begin{tabular}{|c|c|c|c|c|c|c|}
    \hline
        \multicolumn{2}{|c|}{$\Delta S$ }& $\pi/100$ & $\pi/64$ & $\pi/36$ & $\pi/16$ & Order \\ \hline
        \multirow{2}*{Uniform} & g=0 & 4.337E-03 & 5.560E-03 & 8.672E-03 & 1.629E-02 & 0.73 \\ \cline{2-7}
        ~ & g=0.9 & 4.354E-03 & 5.580E-03 & 8.707E-03 & 2.075E-02 & 0.73 \\ \hline
        \multicolumn{2}{|c|}{$\Delta S$ }& $\pi/60$ & $\pi/40$ & $\pi/24$ & $\pi/12$ & Order \\ \hline
        \multirow{2}*{Gaussian} & g=0 & 3.966E-03 & 5.496E-03 & 8.862E-03 & 1.477E-02 & 0.82 \\ \cline{2-7}
        ~ & g=0.9 & 3.978E-03 & 5.518E-03 & 8.902E-03 & 1.491E-02 & 0.82 \\ \hline
    \end{tabular}
\end{table}

\section{Random ordinate method}

In this section, we propose the Random Ordinate Method (ROM) as a potential approach that can achieve a higher order of convergence when the solution has low regularity in the velocity variable, and we demonstrate the possible advantages of this method.

\subsection{Algorithmic details of ROM}

The ROM is based on DOM, but the ordinates are chosen randomly. More precisely, the ROM is performed as follows:
\begin{enumerate}
\item The velocity space $S$ is divided into $n$ cells and each cell is denoted by $ S_{\ell}$ ($\ell=1,\cdots, n$). The maximum area of all $S_{\ell}$ is denoted by $\Delta S=\max_{\ell=1,\cdots,n}|S_{\ell}|$. For example in 1D, $S=[-1,1]$, if uniform mesh is employed, then $\Delta S=2/n$.
\item Sample randomly one ordinate from each cell with uniform probability. Denote $\mathbb{V}^\xi=\{\b{u}_{1}, \cdots, \b{u}_{n}\}$ as the tuple of random ordinates and $V^{\xi}$ as the index set.
\item Solve the resulting discrete ordinate system with the randomly chosen velocity directions.
\begin{equation}\label{equ_random}
\b{u_{\ell}} \cdot \nabla \psi_{\ell}(\b{z})+\sigma_{T}(\b{z}) \psi_{\ell}(\b{z})=\sigma_{S}(\b{z}) \sum_{\ell'=1}^{n}\omega_{\ell'} P_{\ell',\ell}\psi_{\ell'}(\b{z}) +q_{\ell}(\b{z}),\quad \b{u_{\ell}}\in \mathbb{V}^{\xi},
\end{equation}
subject to the boundary conditions
\begin{equation}\label{equ_random_bd}
	\psi_{\ell}(\boldsymbol{z})=\psi_{\Gamma}^{-}(\boldsymbol{z}, \boldsymbol{u_{\ell}}), \quad \boldsymbol{z} \in \Gamma^{-}=\partial \Omega, \quad \boldsymbol{u_{\ell}} \in \mathbb{V}^{\xi}, \quad \boldsymbol{u_{\ell}} \cdot \boldsymbol{n}_{\boldsymbol{z}}<0.
\end{equation}
\end{enumerate}
Since $\boldsymbol{u_{\ell}}$ are now randomly chosen, to guarantee the solution accuracy, one has to determine the corresponding discrete weights $\omega_\ell$ and discrete scattering kernel $P_{\ell',\ell}$. We will simply choose the following approximation:
\begin{gather}
\frac{1}{|S|}\int_{S} P(\b{u'}, \b{u}_\ell)\psi(\b{z}, \b{u'})\,d\b{u'}\approx \sum_{\ell'=1}^{n}\omega_{\ell'} P_{\ell',\ell}\psi_{\ell'}(\b{z}),
\end{gather}
where
\begin{equation}\label{eq:wl}
\omega_{\ell'} =\frac{|S_{\ell'}|}{|S|},
\qquad P_{\ell', \ell}= P(\b{u}_{\ell'},\b{u}_\ell).
\end{equation}
Such a choice clearly satisfies $\sum_{\ell=1}^{n}\omega_{\ell}=1$.
If $\ell'\neq \ell$, then
\[
\E_{\mu_{\ell'}} P_{\ell',\ell}\psi_{\ell'}(\b{z})
=\frac{1}{|S_{\ell'}|}\int_{S_{\ell'}}
P(\b{u}', \b{u}_{\ell})\psi(\b{z}, \b{u}')d\b{u}'.
\]
Therefore, the expectation of the weighted summation provides a good approximation to the integration on the right hand side of \eqref{eq:equ_1}. More precisely,
$$\begin{aligned}
\sum_{\ell'=1}^n\omega_{\ell'}\E_{\mu_{\ell'}} P_{\ell',\ell}\psi_{\ell'}(\b{z})
=&\frac{1}{|S|}\int_{S} P(\b{u'}, \b{u}_\ell)\psi(\b{z}, \b{u'})\,d\b{u'}\\
&-\frac{|S_\ell|}{|S|}\int_{S_\ell} P(\b{u'}, \b{u}_\ell)\psi(\b{z}, \b{u'})\,d\b{u'}+\omega_\ell P(\b{u}_\ell, \b{u}_\ell)\psi(\b{z}, \b{u}_\ell)\\
=&\frac{1}{|S|}\int_{S} P(\b{u'}, \b{u}_\ell)\psi(\b{z}, \b{u'})\,d\b{u'}+O(|S_\ell|^2D(S_{\ell})).
\end{aligned}$$
where $D(S_{\ell})$ means the diameter of the $\ell$-th cell. By employing such a straightforward strategy, there exists the possibility that $\sum_{\ell'}P_{\ell, \ell'}\omega_{\ell'}\neq 1$, potentially leading to a violation of mass conservation at the discrete level. Nonetheless, the positivity of the solution is maintained. Since we consider only $O(1)$ total and scattering cross sections in this current work, only a small error is introduced.

According to \cite{Carlson1955SOLUTION,Edward1989Asymptotic}, symmetric ordinates perform better, especially when there are multiscale parameters in the computational domain. Although we do not consider multiscale parameters in this current paper, symmetric ordinates are used in the ROM. More precisely, in slab geometry, let $n=2m$ and $S_1\cup S_2\cup\cdots\cup S_m=[-1,0]$. $\mu_{\ell}$ ($\ell=1,\cdots,m$) are randomly sampled from $S_{\ell}$ with a uniform distribution. Then
\[\mu_{m+\ell}=-\mu_{m+1-\ell},\quad \mbox{for $\ell=1,2,\cdots,m$}.\]
In the X-Y geometry, let $n=4m$. $S_1\cup S_2\cup\cdots\cup S_{m}$ is the 1/4 disk in the first quadrant, and $\mu_{\ell}$ ($\ell=1,\cdots,m$) are randomly sampled from $S_{\ell}$ with a uniform distribution. The symmetric ordinates indicate that, for $\ell=1,\cdots,m$,
\begin{subequations}\label{eq:2Dordinates_symmetry}
\begin{gather}
\zeta_{\ell}=\zeta_{\ell+m}=\zeta_{\ell+2m}=\zeta_{\ell+3m},\\
\theta_{\ell}=\pi-\theta_{\ell+m}=\pi+\theta_{\ell+2m}=2\pi-\theta_{\ell+3m}.
\end{gather}
\end{subequations}
\begin{remark}
The selection of $P_{\ell', \ell}$ in equation \eqref{eq:wl} may not satisfy the condition $\sum_{\ell'} P_{\ell,\ell'} \omega_{\ell'} = 1$ for anisotropic cases. This condition is crucial for mass conservation in time-evolutionary problems and within the diffusion regime (when $\frac{\sigma_S}{\sigma_T} \approx 1$) \cite{Edward1989Asymptotic}, but it is less critical for the steady state problem when $\frac{\sigma_S}{\sigma_T}$ deviates from $1$. One could normalize $P_{\ell', \ell}$ to ensure this property is strictly maintained, and the impact of such a modification is expected to be minimal. At the current stage, we use \eqref{eq:wl}. The necessary modifications to maintain mass conservation and their corresponding effects will be investigated in our future work. 
\end{remark}

\subsection{Formal analysis of the bias and mean error}\label{Sec:The convergence of ROM}
In this subsection, we perform a formal analysis of ROM, leaving the rigorous analysis for our subsequent works. For a given quadrature $\mathbb{V}^\xi$, one can measure the numerical errors by the difference between $\phi^\xi(\b{z}) = \sum_{\b{u}_\ell \in \mathbb{V}^\xi} \omega_\ell \psi_\ell(\b{z})$ and the reference average density $\phi(\b{z}) = \fint_S \psi(\b{z}, \b{u})$. As demonstrated in Section \ref{Motivation}, when the solution regularity is low in velocity space, the convergence orders of the given quadratures are also low. At first glance, it may seem that there is no benefit to using the ROM.

As in \cite{NEURIPS2019_eb86d510}, we use a 1D integration to illustrate the advantage. Consider the following integration operator:
\begin{equation}
\mathcal{I}(f) = \int_{s_L}^{s_R} f(s') \, ds'.
\end{equation}
Suppose that the function $f(s)$ defined on the interval $[s_L, s_R]$ is a Lipschitz continuous function with a Lipschitz constant $L$. Let $n$ be an integer, and let $s_1 = s_L$, $s_{n+1} = s_R$. Consider a fixed mesh with $n+1$ grid points $\{s_i\}_{i=1}^{n+1}$ whose maximum mesh size is $h$. The integration $\mathcal{I}(f)$ can be approximated by
\[
\mathcal{I}_n(f) = \sum_{i=1}^n f(\xi_i)(s_{i+1} - s_i), \quad \xi_i \in [s_i, s_{i+1}], \quad i = 1, \ldots, n.
\]
The error of this approximation is given by
\[
\sup_f |\mathcal{I}_n(f) - \mathcal{I}(f)| \sim h.
\]
If we consider the randomized quadrature where the sample point $\xi_i$ is uniformly distributed inside $[s_i, s_{i+1}]$, then $\mathbb{E}(\mathcal{I}_n^\xi(f)) = \mathcal{I}(f)$ and the expected value of the distance between $\mathcal{I}_n^\xi(f)$ and $\mathcal{I}(f)$, referred to as the "error," is given by
\[
\sup_f \mathbb{E} |\mathcal{I}_n^\xi(f) - \mathcal{I}(f)| \leq \sup_f \sqrt{\mathrm{Var}(\mathcal{I}_n^\xi(f))} \sim \sqrt{\frac{h^2}{n}}.
\]
Here, the order $\sqrt{\frac{h^2}{n}}$ is due to the fact that the integrals on different subintervals are independent, and the variance for the integral on each subinterval is comparable to the square of the diameter of the interval \cite{ErroranalysisSTENGL1995}.
If we take a uniform grid such that $h = \frac{1}{n}$, this error reduces to $\frac{1}{n^{\frac{3}{2}}}$.
This calculation implies that for a randomized algorithm, the typical order is in fact $n^{-3/2}$, which is better compared to the deterministic one, which is first order for Lipschitz functions (not necessarily smooth). The reason is simple: for any fixed quadrature, there is a choice of functions such that the convergence is only first order. However, for a fixed function, if we consider all possible quadratures, most of the quadratures can achieve an order close to $\frac{3}{2}$.

For ROM, we use similar procedure as described above, where we replaced the exact integral with the randomized integral. The difference is that for the transport equation, $\psi$ is not given in advance; instead, it is solved using a numerical method. Consequently, such an approximation could introduce additional error. Due to this reason, we introduce two quantities to characterize the errors.
The first quantity is the bias:
\begin{equation}
\mathcal{B} := \left\| \mathbb{E}[\phi^{\xi}] - \phi \right\|,
\end{equation}
which gives the distance between the expected value of $\phi^{\xi}$ and the reference solution. Here, the norm is a chosen suitable norm for the functions $\phi: [x_L, x_R] \to \mathbb{R}$. This quantity characterizes the systematic error of the method, which will not vanish if we simply increase the number of experiments. The second quantity is the expected single-run error:
\begin{equation}
\mathcal{E} := \mathbb{E} \left[ \left\| \phi^{\xi} - \phi \right\| \right].
\end{equation}
Note that this expected single-run error could be used to control the variance of the method by Hölder's inequality.

Below, we rewrite the transport equation into the following abstract form:
\begin{equation}
\mathcal{L}\psi = \mathcal{I}(\psi) + Q,
\end{equation}
where $\mathcal{L}$ is the differential operator on the left-hand side of \eqref{eq:equ_1}, $\mathcal{I}(\cdot)$ is the integral operator on the right-hand side, and $Q$ represents the contributions of the source term and the boundary values. The ROM can be written abstractly as
\[
\mathcal{L}\psi^{\xi} = \mathcal{I}_n^{\xi}(\psi^{\xi}) + Q.
\]
Suppose that $\mathcal{I}$ is small compared to $\mathcal{L}$, in the sense that we may expand the solution operator formally as
\[
(\mathcal{L} - \mathcal{I})^{-1} = \mathcal{L}^{-1} + \mathcal{L}^{-1}\mathcal{I}\mathcal{L}^{-1} + \mathcal{L}^{-1}\mathcal{I}\mathcal{L}^{-1}\mathcal{I}\mathcal{L}^{-1} + O(\|\mathcal{I}\mathcal{L}^{-1}\|^3),
\]
so that
\begin{equation}\label{eq:psiexpansion}
\psi = \mathcal{L}^{-1}Q + \mathcal{L}^{-1}\mathcal{I}\mathcal{L}^{-1}Q + \mathcal{L}^{-1}\mathcal{I}\mathcal{L}^{-1}\mathcal{I}\mathcal{L}^{-1}Q + O(\|\mathcal{I}\mathcal{L}^{-1}\|^3).
\end{equation}
This expansion corresponds to the classical source iteration method \cite{Larsen2010}.
A similar expansion is valid for $(\mathcal{L} - \mathcal{I}_n^{\xi})^{-1}$. Then, one finds that
\begin{equation}\label{eq:psixiexpansion}
\psi^{\xi} = \mathcal{L}^{-1}Q + \mathcal{L}^{-1}\mathcal{I}_n^{\xi}\mathcal{L}^{-1}Q + \mathcal{L}^{-1}\mathcal{I}_n^{\xi}\mathcal{L}^{-1}\mathcal{I}_n^{\xi}\mathcal{L}^{-1}Q + O(\|\mathcal{I}_n^{\xi}\mathcal{L}^{-1}\|^3),
\end{equation}
\begin{equation}\label{eq:expectpsiexpansion}
\mathbb{E} \psi^{\xi} = \mathcal{L}^{-1}Q + \mathcal{L}^{-1}\mathbb{E}[\mathcal{I}_n^{\xi}]\mathcal{L}^{-1}Q + \mathbb{E}[\mathcal{L}^{-1}\mathcal{I}_n^{\xi}\mathcal{L}^{-1}\mathcal{I}_n^{\xi}\mathcal{L}^{-1}Q] + O(\|\mathcal{I}_n^{\xi}\mathcal{L}^{-1}\|^3).
\end{equation}
Comparing \eqref{eq:psiexpansion} and \eqref{eq:expectpsiexpansion}, since $Q$ is given and $\mathcal{L}^{-1}$ is a linear operator,
\[
\mathcal{L}^{-1}\mathbb{E}[\mathcal{I}_n^{\xi}]\mathcal{L}^{-1}Q = \mathcal{L}^{-1}\mathcal{I}\mathcal{L}^{-1}Q,
\]
while
\[
\mathbb{E}[\mathcal{L}^{-1}\mathcal{I}_n^{\xi}\mathcal{L}^{-1}\mathcal{I}_n^{\xi}\mathcal{L}^{-1}Q] \neq \mathcal{L}^{-1}\mathcal{I}\mathcal{L}^{-1}\mathcal{I}\mathcal{L}^{-1} Q.
\]
Note that since $\phi$ is the integral of $\psi$, the bias and error for $\phi$ are comparable to those for $\psi$. The difference between $\psi$ and $\mathbb{E} \psi^{\xi}$ gives the leading order of the bias. Clearly, the bias is of the same order as the variance of the randomized integral operator $\mathcal{I}_n^{\xi}$:
\[
\mathcal{B} \lesssim \sup_f \mathrm{Var}(\mathcal{I}_n^{\xi}(f)),
\]
which would be $\frac{h^2}{n}$ for the slab geometry. Comparing \eqref{eq:psiexpansion} and \eqref{eq:psixiexpansion}, for the expected single-run error, one has
\[
\mathcal{E} \sim \|\mathcal{L}^{-1}(\mathcal{I}_n^{\xi} - \mathcal{I})\mathcal{L}^{-1}Q\| \lesssim \sup_f \sqrt{\mathrm{Var}(\mathcal{I}_n^{\xi}(f))}.
\]
Hence, for the slab geometry, we expect that ROM can achieve a convergence order of $\frac{3}{2}$, and the expectation of multiple runs gives a third-order convergence.

One can generalize the analysis to higher dimensions.  If one divides the region $S$ into $n$ cells, the random ordinates are independent and the variance for the random quadrature $\mathcal{I}_n^{\xi}$ is like $\max_{\ell}D(S_{\ell})^2/n$, where $D(S_{\ell})$
is the diameter of the $\ell$-th cell. Consequently, 
the bias scales like 
\[
\|\phi-\E\phi^{\xi}\|
\sim \frac{\max_{\ell}D(S_{\ell})^2}{n}\sim \frac{1}{n^{1+2/d}}.
\]
The error scales like 
\[
\E\|\phi-\phi^{\xi}\|\sim \sqrt{\frac{\max_{\ell}D(S_{\ell})^2}{n}}
\sim \frac{1}{n^{1/2+1/d}}.
\]
For $d=2$, a typical order of error for ROM is then $1$, and the bias scales like $O(n^{-2})$ so the order is $2$. The order of bias also improves compared to DOM. Hence, ROM could improve accuracy.

\subsection{Comparison of the computational cost of DOM and ROM}

We see from the formal analysis above that the convergence order of bias is higher even if the solution regularity in velocity space is low. Therefore, if we take more samples of $\mathbb{V}^\xi$, run the system \eqref{equ_random} multiple times in parallel, and then take the expectation $\mathbb{E}\phi^{\xi}$, the solution accuracy can be improved.

In the X-Y geometry (i.e., $d=2$), if the required accuracy for the density $\phi(x,y)$ is $\epsilon$, and if the convergence order of DOM is $1$, the mesh size needs to be $O(\epsilon)$. Consequently, the required number of ordinates is $n = O(\epsilon^{-1})$. The standard source iteration method with an anisotropic scattering kernel requires, in each iteration step, the multiplication of the discrete anisotropic scattering kernel, which is an $n \times n$ matrix, with the density flux of the current iteration step, which is an $n \times 1$ vector. This results in $n^2$ multiplications. Therefore, one can assume that the computational cost of solving an $n$ ordinates DOM system is $O(n^2) = O(\epsilon^{-2})$.
If one solves the large coupled linear system using other iterative solvers, the computational complexity can be higher than $O(\epsilon^{-2})$.

On the other hand, the bias $\mathcal{B}$ of ROM is of order $2$. If we repeat the simulation for $t$ times, the error of the Monte Carlo approximation is 
$O(\sqrt{\mathrm{Var}/t+\mathcal{B}^2})$. The variance can be controlled by the square of the error, which is of order $1$. Hence, the number of ordinates for each run could be $O(\epsilon^{-\frac{1}{2}})$ so that $\mathcal{B}=O(\e)$ and $\mathrm{Var}=O(\e)$. Then, taking $t=O(\e^{-1})$ will make the average error $O(\e)$. Since the complexity for each run is $O([\epsilon^{-\frac{1}{2}}]^2)=O(\e^{-1})$, the total complexity is thus $O(t\e^{-1})=O(\e^{-2})$. The total computational costs of DOM and ROM are comparable.

If one chooses to use $O(\e^{-1})$ ordinates for DOM, then a typical single run could yield $O(\e)$ error, but it suffers from the ray effect. On the other hand, choosing $O(\e^{-1/2})$ ordinates and repeat $t=O(\e^{-1})$ times could result in the same error, but the ray effect could be improved since the total number of velocity directions is $O(\e^{-3/2})$. Moreover, ROM is easy to parallel.

 The key takeaway of this section is that when the solution exhibits low regularity and the DOM convergence rate is at most first-order, the ROM demonstrates its superiority. This theoretical conclusion aligns with the numerical tests in the X–Y geometry (Figure \ref{fig:costvsbias}), where the ROM yields smaller errors at equivalent computational costs. For the slab-geometry results (Figure \ref{fig:1dcostvsbias}): When the DOM converges faster than first order, the ROM provides negligible benefit; by contrast, When the DOM convergence order approaches 1, the ROM exhibits a clear advantage.

\section{Numerical experiments}
The numerical performance of the ROM is presented in this section. The numerical convergence orders of errors and bias in both slab and X-Y geometries are displayed, while its ability to mitigate the ray effect is demonstrated through a lattice problem.

\subsection{The slab geometry}

In this subsection, we consider following example that when the inflow boundary conditions have low regularity in velocity space.
Though there is no ray effect in slab geometry, we will show the convergence orders of DOM are low verify the error estimates of ROM in Section \ref{Sec:The convergence of ROM} .

\begin{example}
\label{example:1Dthreecases}
Let the computational domain, total and absorption cross sections, scattering kernel and source term be respectively
$$
     x\in [0,1],
     \quad
	\sigma_{T}(x)=10x^2+1,
	\quad 
	\sigma_{S}(x)=5x^2+0.5,\quad P(\mu',\mu)=1,
	\quad q(x)=1+x.
$$
We consider three different inflow boundary conditions:
\begin{enumerate}
 \item[Case 1:]   Continuous and smooth inflow boundary conditions:
 \begin{equation}\label{eq:2.12}
     \psi(0,\mu)=3 \mu ,  \quad \mu>0;
     \qquad
     \psi(1,\mu)=-5 \mu,\quad \mu<0.
 \end{equation}
 \item[Case 2:]   Continuous but non-differentiable inflow boundary conditions:
 \begin{equation}
  \psi(0,\mu)=
 \begin{cases}
    \frac{4}{3}-\mu, & \frac{1}{3}<\mu<1, \\
     3 \mu,  & 0<\mu \le \frac{1}{3} ,
 \end{cases}
\qquad  
    \psi(1,\mu)=
 \begin{cases}
    2+\mu, & -1 <\mu<-\frac{1}{3}, \\
    -5 \mu,  & -\frac{1}{3} \le \mu < 0.
 \end{cases}\label{eq:2.13}
 \end{equation}
 \item[Case 3:] Discontinuous inflow boundary conditions:
 \begin{equation}\label{eq:2.14}
  \psi(0,\mu)=
 \begin{cases}
    3-\mu, & \frac{1}{3}<\mu<1, \\
     3 \mu,  & 0<\mu \le \frac{1}{3}, 
 \end{cases}
\qquad  
    \psi(1,\mu)=
 \begin{cases}
    4+\mu, & -1 <\mu<-\frac{1}{3}, \\
    -5 \mu,  & -\frac{1}{3} \le \mu < 0.
 \end{cases}
 \end{equation}
\end{enumerate}
\end{example}

We employ second-order finite difference spatial discretization in \cite{ShiA} to obtain the numerical results.
The number of spatial cells is fixed to be $I=50$, and the grid points are $x_i$ ($i=0,1,\cdots,I$). Let $\psi_{\ell}(x_i)$ be the solution to the equation, and the average density $\phi(x_i)=\sum_{{\ell}\in V}\omega_{\ell} \psi_{\ell}(x_i)$. The reference solution is computed using $2560$ ordinates. The $\ell^2$ errors of the numerical solutions with different ordinates are defined by
\begin{equation*}
\mathcal{E}_{DOM}=\sqrt{\frac{1}{I+1}\sum_{i=0}^{I}\mid\phi(x_i)-\phi^{ref}(x_i)\mid^{2}}.
\end{equation*}

\begin{figure}[htbp]
	\centering
\subfigure[Uniform quadrature]
{
		\includegraphics[width=0.45\linewidth]{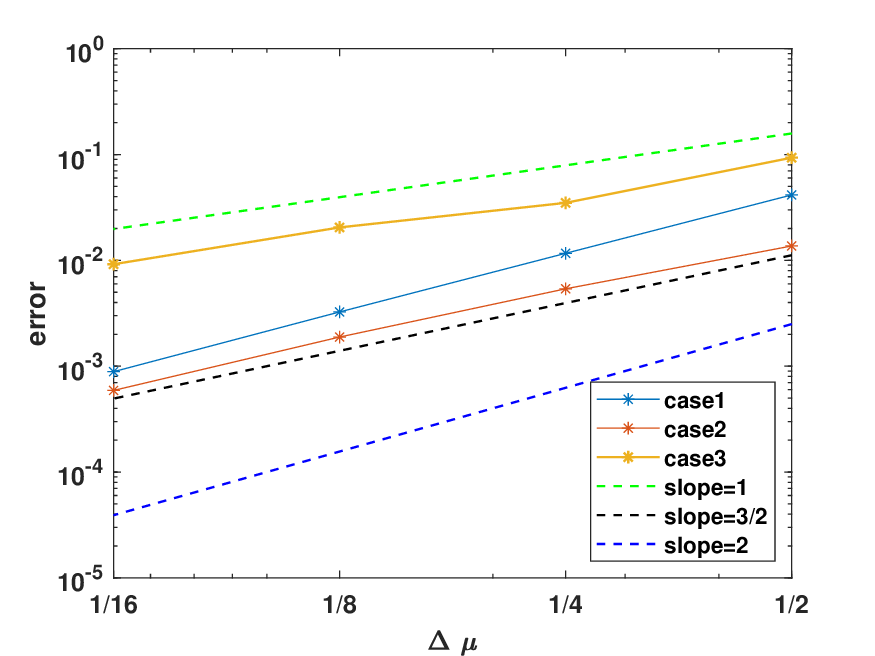}
}
\subfigure[Gaussian quadrature]
{
		\includegraphics[width=0.45\linewidth]{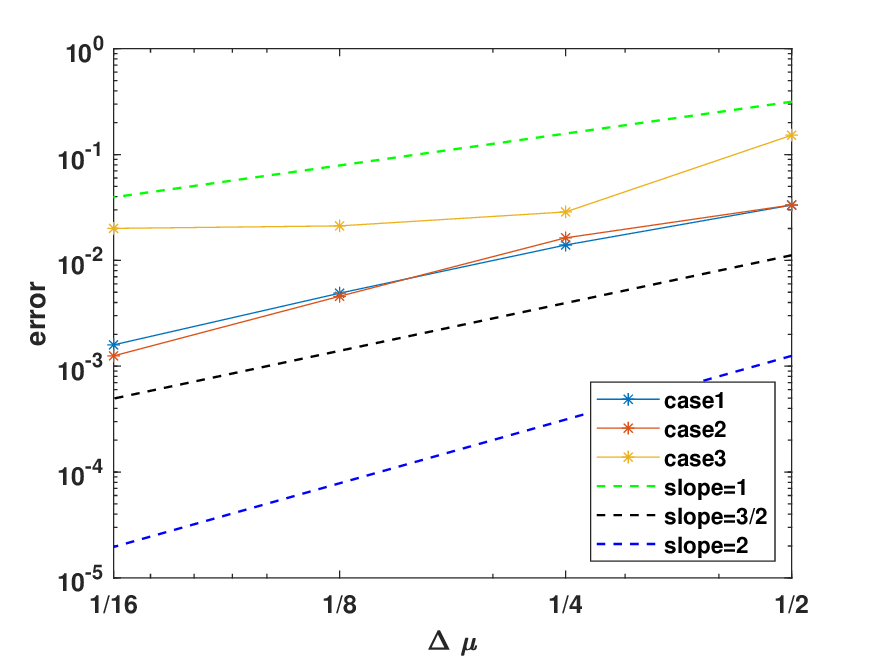}
}

		\caption{Example \ref{example:1Dthreecases}: Convergence orders of DOM with different inflow boundary conditions in \eqref{eq:2.12}-\eqref{eq:2.14}. (a): Uniform quadratures of different sizes;   (b): Gaussian quadratures of different sizes. Here $\Delta\mu=\frac{1}{M}$. }
		\label{fig:section2_1D_DOMconvergence}
\end{figure}

The convergence orders of DOM for different cases are shown in Figure \ref{fig:section2_1D_DOMconvergence}. One can observe that for both Uniform and Gaussian quadratures, the convergence orders decrease from 2 to 1 when the inflow boundary conditions change from Case 1 to Case 3. Therefore, one cannot expect a high convergence order when the regularity of the inflow boundary conditions is low. In particular, Gaussian quadrature does not reach spectral convergence as shown in Figure \ref{fig:section2_1D_DOMconvergence}(b). This is because, though the inflow boundary condition is smooth in $\mu$, at the boundary, the solution jumps at $\mu=0$. Gaussian quadrature does not provide spectral convergence for solutions with a jump at $\mu=0$. 

To apply ROM in the slab geometry, we equally divide the velocity interval $[-1,1]$ into $n$ cells. The $\ell$th cell is denoted by $S_{\ell}=[-1+\frac{2\ell-2}{n},-1+\frac{2\ell}{n}]$, where $\ell=1,2,\cdots,n$. We consider only even $n$, and when $S_{\ell}\subset [-1,0]$, one ordinate $\mu_{\ell}$ is sampled randomly from the cell $S_{\ell}$ with uniform probability. When $S_{\ell}\subset[0,1]$, $\mu_{\ell}=-\mu_{n+1-\ell}$. The weights are $\omega_{\ell}=2/n$ for all $\mu_{\ell}$. 

We define the error between the reference solution and the numerical solutions of ROM with $I$ spatial cells by
\begin{equation}\label{eq:norm-error}
	\mathcal{E}=\mathbb{E}\parallel \phi^{\xi}(x)-\phi^{ref}(x)\parallel_2=\mathbb{E}\left(\frac{1}{I+1}\sum_{i=0}^{I}\mid\phi^{\xi}(x_i)-\phi^{ref}(x_i)\mid^{2}\right)^{\frac{1}{2}},
\end{equation}
where $\phi^\xi(x_i)=\sum_{m=1}^{n}\omega_m \psi^\xi_{m}(x_i)$.
The bias of ROM is defined by
\begin{equation}\label{eq:norm-bias}
	 \mathcal{B}=\parallel \mathbb{E}\phi^{\xi}(x)-\phi^{ref}(x)\parallel_2=\left(\frac{1}{I+1}\sum_{i=0}^{I}\mid\mathbb{E}\phi^{\xi}(x_i)-\phi^{ref}(x_i)\mid^{2}\right)^{\frac{1}{2}}.
\end{equation}

\begin{figure}[htp]
	\centering
\subfigure[Example \ref{example:1Dthreecases}]
{
		\includegraphics[width=0.45\linewidth]{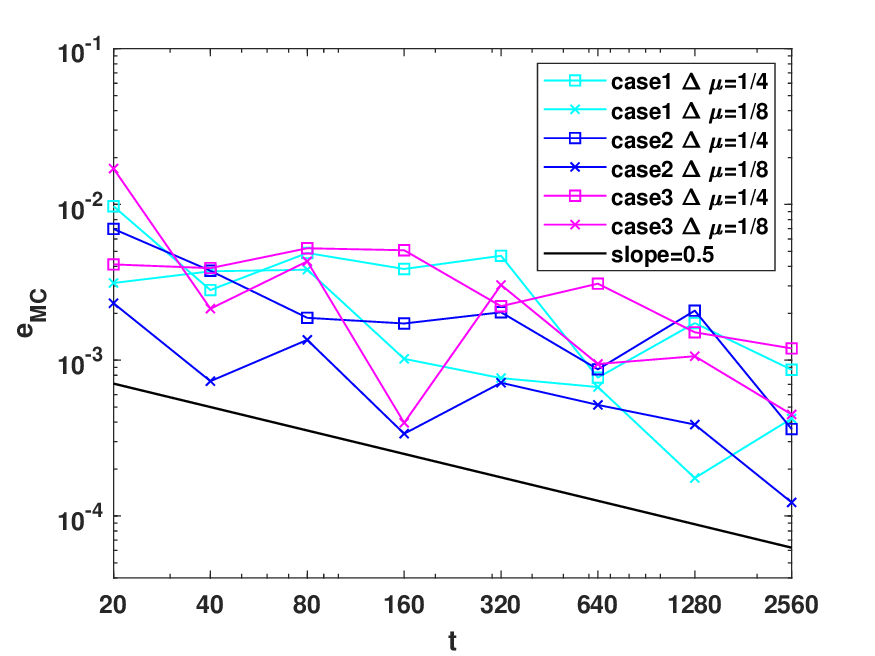}
        \label{fig:1dMC}
}
\subfigure[ Example \ref{2D disk_ray effect}]
{
		\includegraphics[width=0.45\linewidth]{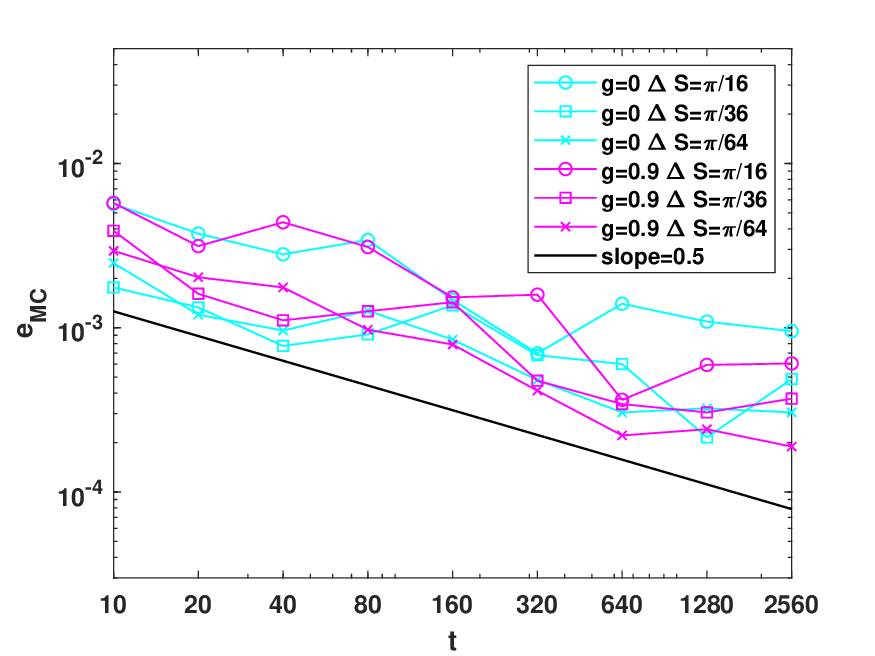}
        \label{fig:2dMC}
}
		\caption{ The convergence of ROM with respect to the number of sampled simulations $t$. (a): the MC error defined in \eqref{eq:1dMCerror} for different cases in Example \ref{example:1Dthreecases}; (b): the MC error for different scattering kernel in Example \ref{2D disk_ray effect}. $t$ is the number of sampled simulations in ROM.}
		\label{fig:MC}
\end{figure}

\paragraph{Convergence with respect to the number of sampled simulations}
In ROM, let $\phi^{\xi_s}(x)$ represent the total density for the $s$-th sample. The average of $\phi^{\xi_s}(x)$ over $t$ sampled simulations is given by
\begin{equation}
    \overline{\phi}^{\xi}(x) = \frac{1}{t} \sum_{s=1}^t \phi^{\xi_s}(x).
\end{equation}
The Monte Carlo (MC) error is then expressed as
\begin{equation}\label{eq:1dMCerror}
    e^t_{MC} = \left\| \overline{\phi}^{\xi}(x) - \mathbb{E}[\phi^{\xi}(x)] \right\|_2 = \left( \frac{1}{I+1} \sum_{i=0}^{I} \left| \frac{1}{t} \sum_{s=1}^t \phi^{\xi_s}(x_i) - \mathbb{E}[\phi^{\xi}(x_i)] \right|^2 \right)^{\frac{1}{2}}.
\end{equation}
Here, $\mathbb{E}[\phi^{\xi}(x)]$ is obtained by $\frac{1}{20480} \sum_{s=1}^{20480} \phi^{\xi_s}(x)$.

Figure \ref{fig:MC} presents the convergence orders of ROM with respect to the number of sampled simulations $t$.
Figure \ref{fig:1dMC} corresponds to Example \ref{example:1Dthreecases} in slab geometry, where three different inflow boundary conditions given by equations \eqref{eq:2.12}--\eqref{eq:2.14} with varying regularities are considered.
It can be observed that the convergence order of the MC error with respect to the number of sampled simulations $t$ is 0.5, which is consistent with the behavior of the Monte Carlo method. 

\begin{figure}[htp]
	\centering
\subfigure[error $\mathcal{E}$]
{
		\includegraphics[width=0.45\linewidth]{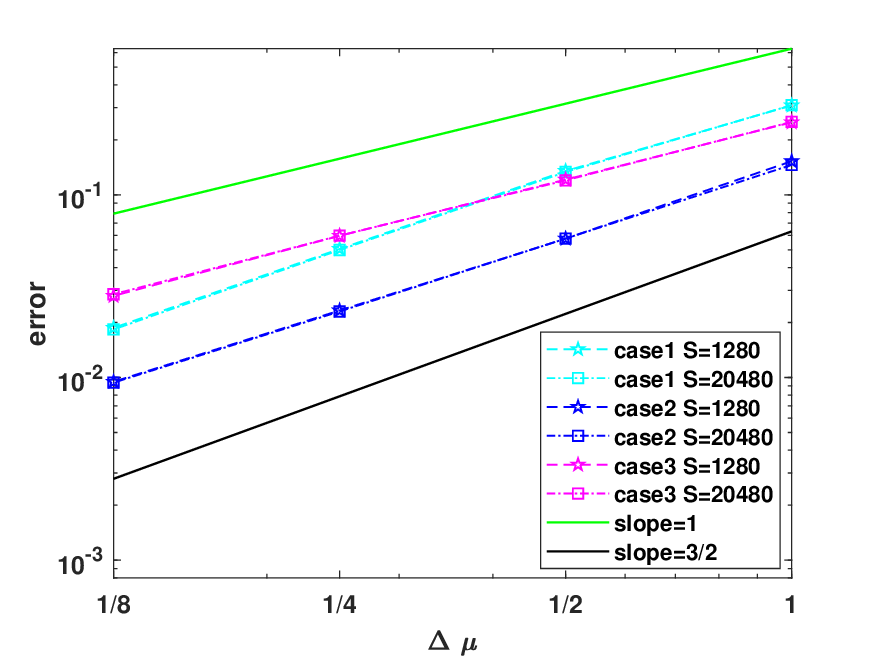}
}
\subfigure[bias $\mathcal{B}$]
{
		\includegraphics[width=0.45\linewidth]{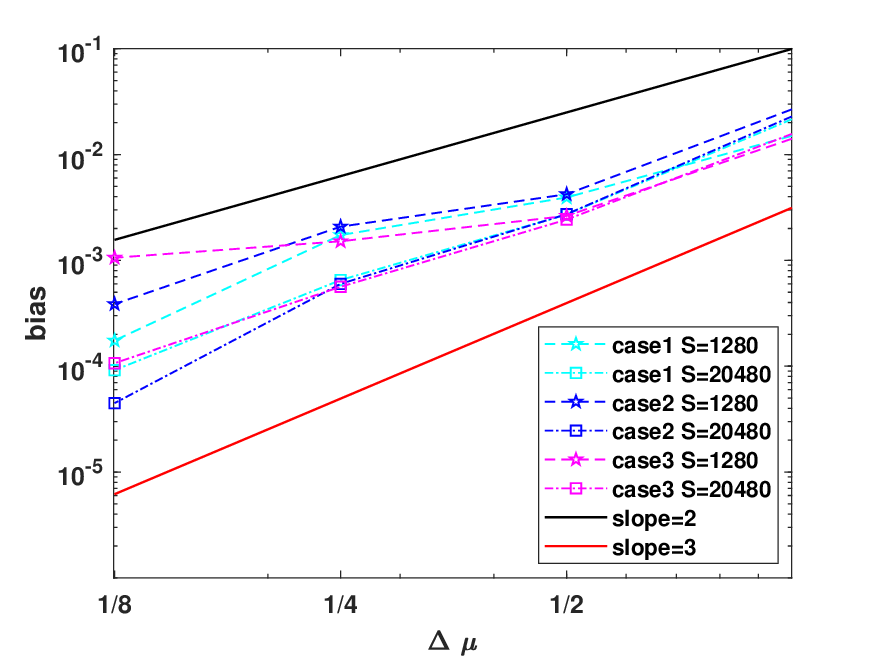}
}
		\caption{ The convergence orders of ROM in velocity in slab geometry. Here $\Delta \mu=\frac{2}{n}$,$n=2,4,8,16$. (a): the errors  defined in \eqref{eq:norm-error} for different cases in \eqref{eq:2.12}-\eqref{eq:2.14}; (b): the MC error for different cases in \eqref{eq:2.12}-\eqref{eq:2.14}. $t$ is the number of sampled simulations in ROM. The values $t = 1280$ and $20480$ were selected. In practice, a magnitude of $10^3$ usually suffices, but here we need $t$ large enough for the MC error to reach the asymptotic, high-order regime so that the observed slope matches the theoretical value.}
		\label{fig:1D_3case}
\end{figure}

\begin{table}[!ht]
    \caption{Error, bias and slope of convergence curves in Figure \ref{fig:1D_3case} for $t=20480$.}
    \label{Table:1dROM_S20480}
    \centering
    \begin{tabular}{|c|c|c|c|c|c|c|}
    \hline
        \multicolumn{2}{|c|}{$\Delta \mu$ }  &  1/8 &  1/4 &  1/2 & 1     &  Order \\ \hline
        \multirow{3}*{error ($\mathcal{E})$} & Case 1 & 1.83E-02 & 4.99E-02 & 1.33E-01 & 3.10E-01 & 1.37  \\ \cline{2-7}
        ~ & Case 2 & 9.37E-03 & 2.30E-02 & 5.76E-02 & 1.46E-01 & 1.32  \\ \cline{2-7}
        ~ & Case 3 & 2.86E-02 & 5.98E-02 & 1.20E-01 & 2.51E-01 & 1.04  \\ \hline
        \multirow{3}*{bias	($\mathcal{B})$} & Case 1 & 9.18E-05 & 6.51E-04 & 2.72E-03 & 2.17E-02 & 2.57  \\ \cline{2-7}
        ~ & Case 2 & 4.48E-05 & 6.00E-04 & 2.73E-03 & 2.31E-02 & 2.92  \\ \cline{2-7}
        ~ & Case 3 & 1.07E-04 & 5.63E-04 & 2.43E-03 & 1.58E-02 & 2.37 \\ \hline
    \end{tabular}

\end{table}

\paragraph{Convergence with respect to the mesh sizes} Figure \ref{fig:1D_3case} displays the convergence orders of the error and bias of ROM with respect to mesh sizes. The results for different numbers of sampled simulations are presented. To determine the convergence order of the error, a smaller number of samples suffices. As observed from Figure \ref{fig:1D_3case}, as the number of samples increases, the convergence order of the bias gradually increases and eventually stabilizes near the theoretical value. For different cases as described in Example \ref{example:1Dthreecases} in Section 2, Table \ref{Table:1dROM_S20480} shows the error and bias using various mesh sizes when an adequate number of samples is employed. In Figure \ref{fig:1D_3case}, to observe the convergence order without the influence of Monte Carlo errors, $t$ is chosen to be sufficiently large.

\paragraph{Comparison of ROM with DOM} When the inflow boundary conditions are smooth and regular as in Case 1 in \eqref{eq:2.12}, the convergence orders of DOM are $2$, while the convergence orders of the error and bias of ROM are respectively $1.37$ and $2.57$. Due to stochastic noise, this can be considered consistent with the analysis.
When the boundary conditions are continuous but nondifferentiable as in Case 2 in \eqref{eq:2.13}, the convergence orders of the DOM decrease to $1.5$, and the convergence orders of the error and bias of ROM remain $1.32$ and $2.92$, respectively. Moreover, if the boundary conditions are discontinuous at some points as in Case 3 in \eqref{eq:2.14}, the convergence order of DOM decreases to $1$. The convergence order of the error of ROM decreases to $1$, while the bias has a convergence order of $2.37$. In summary, the errors of ROM converge no slower than DOM, and the bias converges faster, especially when the solution regularity is low in the velocity variable. The details of the convergence are shown in Table \ref{Table:1dROM_S20480}.

\begin{figure}[ht]
	\centering
\subfigure[case 2]
{
		\includegraphics[width=0.45\linewidth]{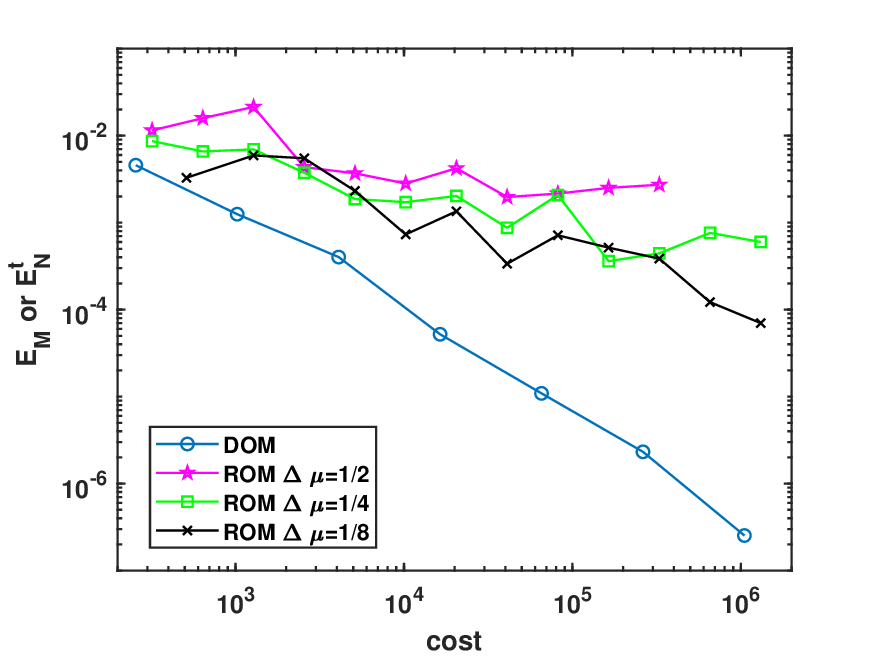}
}
\subfigure[case 3]
{
		\includegraphics[width=0.45\linewidth]{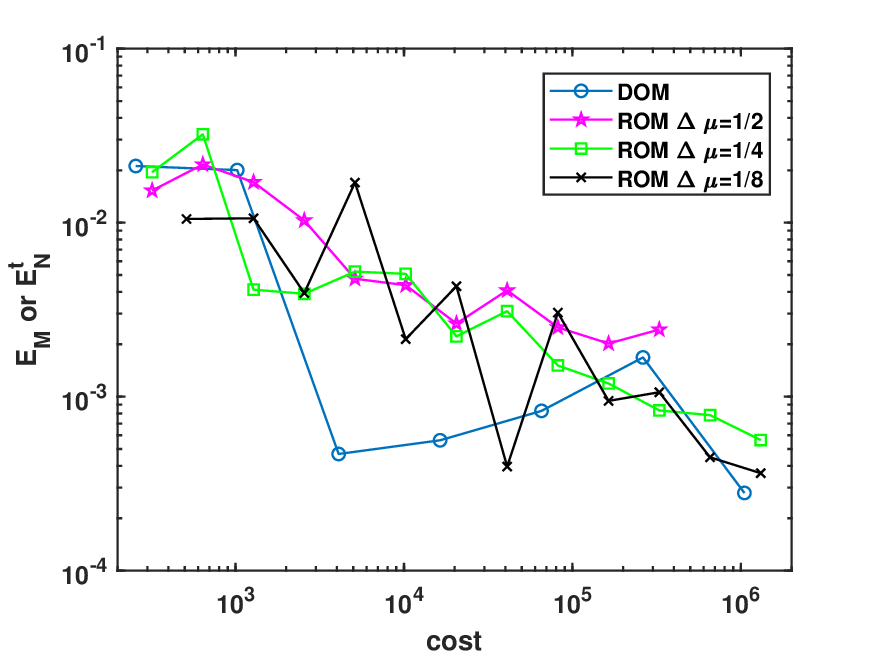}
}

		\caption{Comparison of the DOM and ROM accuracy in terms of the computational cost in slab geometry:
For different numbers of ordinates $M = 2, 4, 8, 16, 32, 64, 128, 256, 512, 1024$, the costs of DOM are assumed to be $M^2$, and the obtained numerical solutions are compared with the reference solution, i.e., the values of $E_M = \|\phi_M - \phi^{ref}\|_2$ are displayed. For ROM, the number of ordinates for each sampled simulation is fixed at $N$, and as the number of sampled simulations $t$ increases from 2 to 20480, the cost is assumed to be $tN^2$. For $N = 4, 8, 16$, the obtained average solution is denoted by $\phi^t_N$, and $E^t_N = \|\phi^t_N - \phi^{ref}\|_2$ are displayed.
(a) Case 2 in Example \ref{example:1Dthreecases}; (b) Case 3 in Example \ref{example:1Dthreecases}. }
		\label{fig:1dcostvsbias}
\end{figure}

\paragraph{Comparison of the complexity.} Assume that the number of ordinates used in DOM and ROM are $M$ and $N$, respectively, and the number of sampled simulations in ROM is $t$. As discussed in Section \ref{Motivation}, when employing classical spatial discretization and the standard source iteration method to solve the discretized system, the cost of DOM with respect to $M$ is $O(M^2)$. Consequently, the cost of each sampled simulation in ROM is $O(N^2)$, and the total cost of ROM is $O(tN^2)$. To compare the cost of DOM and ROM, in DOM, we select different values of $M = 2, 4, 8, 16, 32, 64, 128, 256, 512, 1024$ and compare the obtained numerical solutions with the reference solution obtained with a very fine mesh, i.e., $E_M = \|\phi_M - \phi^{ref}\|_2$. For different $M$, the values of $(M^2, E_M)$ are shown in Figure \ref{fig:1dcostvsbias}. For ROM, we fix $\Delta \mu$ and increase the number of sampled simulations $t$ from $2$ to $20480$. We consider $N = 4, 8, 16$ and the obtained average solution is denoted by $\phi^t_N$. Let $E^t_N = \|\phi^t_N - \phi^{ref}\|_2$, the values of $(tN^2, E^t_N)$ are also displayed in Figure \ref{fig:1dcostvsbias}. We observe that, when the solution regularity is low, the solution accuracy in DOM depends on the locations of the quadrature nodes, but the accuracy of ROM is comparatively stable. In slab geometry, when the solution regularity is low in the velocity variable, the computational costs of DOM and ROM are comparable.

\subsection{The X-Y geometry case}
\begin{figure}[ht]
	\centering
\subfigure[]
{
		\includegraphics[width=0.45\linewidth]{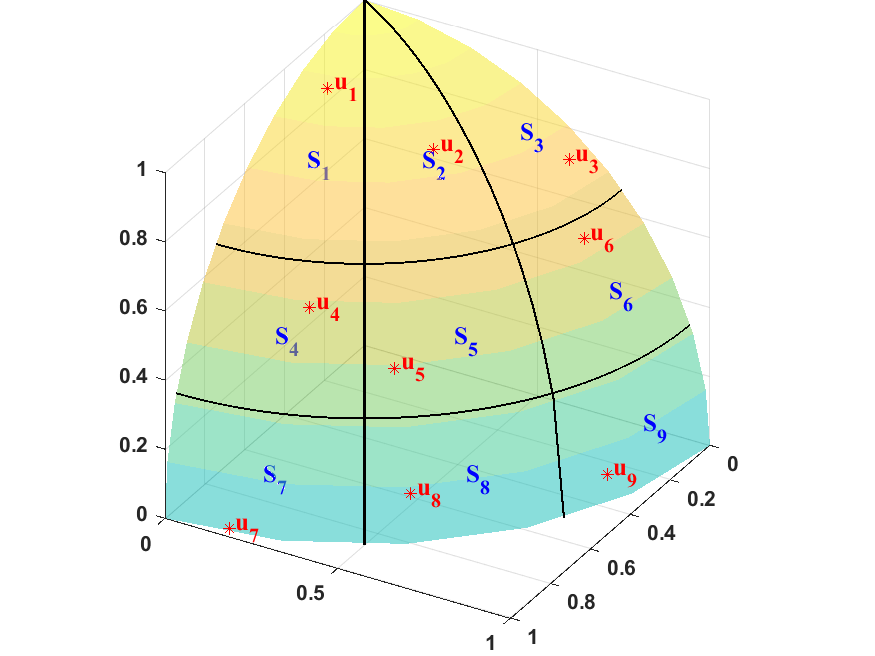}
}
\subfigure[]
{
		\includegraphics[width=0.45\linewidth]{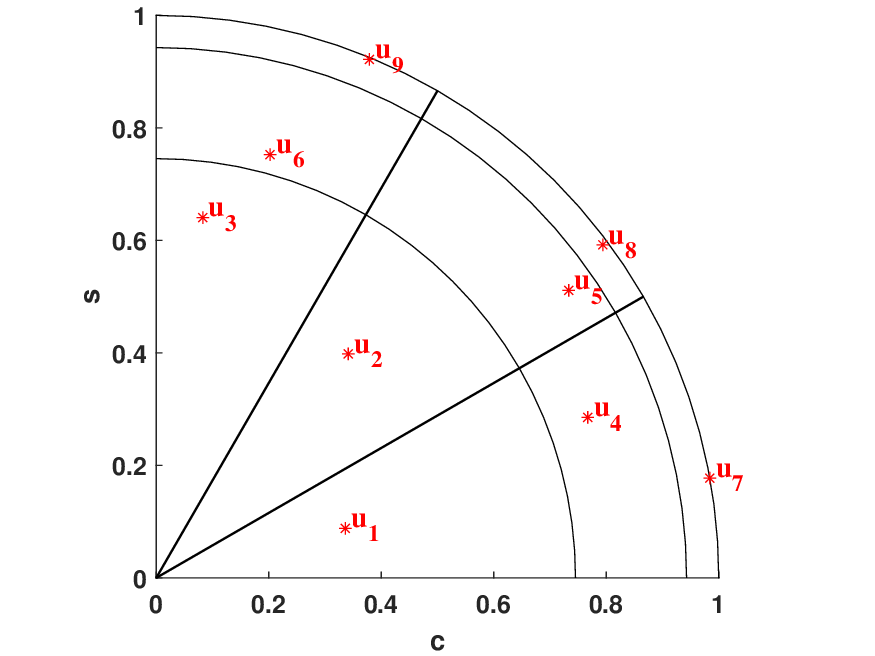}
}

		\caption{Schematic diagram of selected ordinates on the surface of a 3D unit sphere and their corresponding projection to the 2D unit disk of ROM. }
		\label{fig:2dROM_ordinates}
\end{figure}

In X-Y geometry, we show the ordinates in the first quadrant and consider uniform partition in the $(\zeta,\theta)$ plane.    $\zeta\in[0,1]$ and $\theta\in[0,\pi/2]$ are equally divided, the nodes are denoted by $(\zeta_i,\theta_j)=(\frac{N-i+1}{N},\frac{j-1}{2N}\pi)$ for $i,j=1,\cdots,N+1$. Each quadrant has $m=N^2$ cells. For any $\ell=(i-1)N+j$, $i,j=1,\cdots,N$,
$$
S_{\ell}=\{ (\zeta,\theta)| \zeta_{a+1} \leq \zeta \leq \zeta_{a}, \theta_{b} \leq \theta \leq \theta_{b+1}. \}
$$
Randomly sample $\zeta^{\xi}_i \in [\zeta_{i+1},\zeta_{i}]$ and $\theta^{\xi}_j \in [\theta_{j},\theta_{j+1}]$ with uniform probability. The ordinates in other quadrant can be obtained by symmetry as in \eqref{eq:2Dordinates_symmetry}. The weights are chosen to be $\bar\omega_\ell=\frac{1}{4N^2}$.   The ordinates projected to the 2D unit disk are
\begin{equation*}
(c_{\ell},s_{\ell},\bar{\omega}_{\ell})=\Big(\big(1-(\zeta^\xi_i)^2\big)^{\frac{1}{2}}\cos\theta_j^{\xi},\big(1-(\zeta^{\xi}_i)^2\big)^{\frac{1}{2}}\sin\theta_{j}^{\xi},\frac{1}{4N^2}\Big).
\end{equation*}
Let’s take $N = 3$ as an example, the partition and one sample of the chosen discrete ordinates are plotted in Figure \ref{fig:2dROM_ordinates}.
\subsubsection{Convergence order} Similar to section 2.2.2, we use the diamond difference method to discretize the spatial variable. The unknowns at the cell centers are calculated, and the notations are the same as in section 2.2.2. We use the same setup and spatial grid as Example \ref{2D disk_ray effect} in section 2. The numerical error $\mathcal{E}$ and bias $\mathcal{B}$ are similar to \eqref{eq:norm-error} and \eqref{eq:norm-bias}, and the MC error is similar to \eqref{eq:1dMCerror}, expect that the $\ell^2$ norm is given as in \eqref{eq:l2norm2D}.
The reference solution $\phi^{ref}$ is computed by 80400 ordinates by Gaussian quadrature.

Figure \ref{fig:2dMC} displays the convergence orders of ROM with respect to the number of sampled simulations $t$. The numerical results for both isotropic and anisotropic cases, calculated with three different mesh sizes $\Delta S = \pi/16$, $\pi/36$, and $\pi/64$, are shown. The convergence order of the MC error with respect to the number of sampled simulations $t$ is $0.5$.

\begin{figure}[ht]
	\centering
\subfigure[error $\mathcal{E}$]
{
		\includegraphics[width=0.45\linewidth]{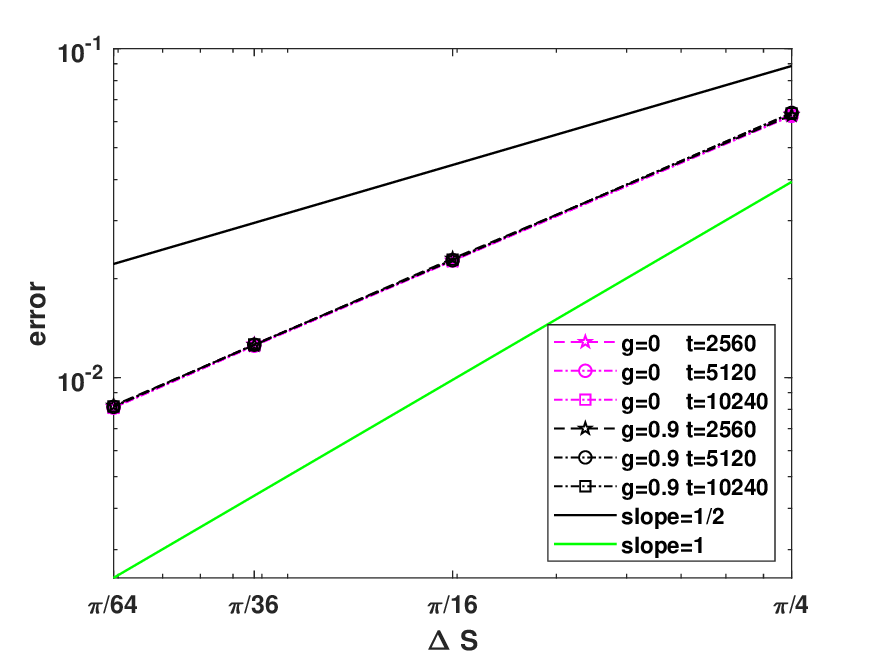}
}
\subfigure[bias $\mathcal{B}$]
{
		\includegraphics[width=0.45\linewidth]{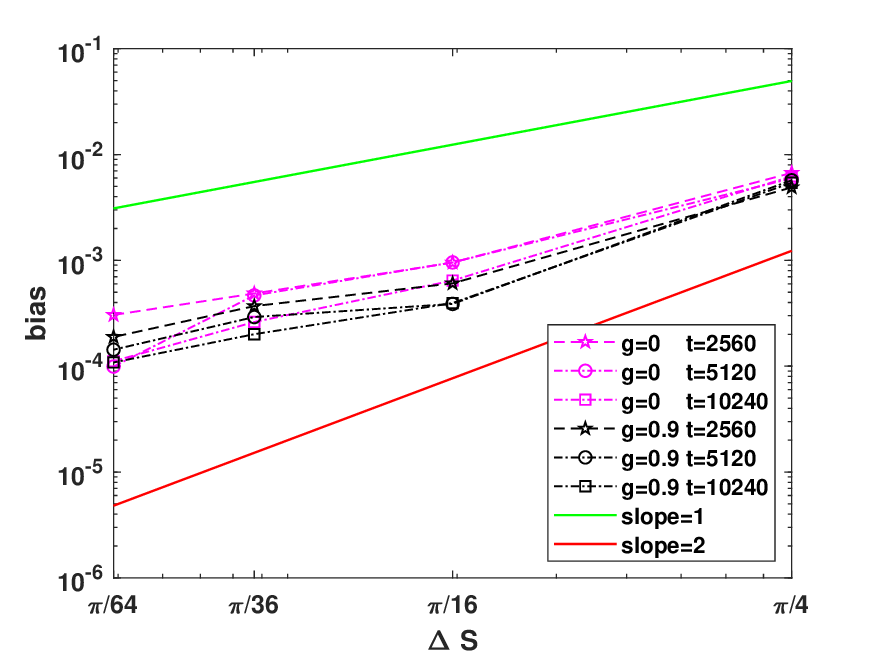}
}

		\caption{The convergence orders of ROM with different velocity meshes in X-Y geometry. Here $\Delta S=\frac{\pi}{4N^2}$,$N=1,2,3,4$.  (a) the errors ($\mathcal{E}$) for isotropic and anisotropic scattering kernel, (b): the bias ($\mathcal{B}$) for isotropic and anisotropic scattering kernel. t is the number of sampled simulations in ROM. The values $t = 2560$ and $10240$ were selected. Here we need $t$ large enough for the MC error to reach the asymptotic, high-order regime so that the observed slope matches the theoretical value.} 
		\label{fig:2D_ROMconvergence}
\end{figure}

 Figure \ref{fig:2D_ROMconvergence} and Table \ref{Table:2dROM_10240} show the convergence order of the error and bias. We can observe that the convergence order of the error is between $0.5$ and $1$, while the convergence order of the bias is between $1$ and $2$ for both isotropic and anisotropic scattering kernels. This is due to the non-smoothness of the solution, which leads to a slight difference between the theoretical value and the actual result.

\begin{table}[!ht]
    \centering
        \caption{Error, bias and slope of convergence curves in Figure \ref{fig:2D_ROMconvergence} for $t=10240$.}
    \label{Table:2dROM_10240}
    \begin{tabular}{|c|c|c|c|c|c|c|}
    \hline
        \multicolumn{2}{|c|}{$\Delta S$ }& $\pi/64$ & $\pi/36$ & $\pi/16$ & $\pi/4$ & Order \\ \hline
         \multirow{2}*{error ($\mathcal{E})$}& iso: g=0 & 8.09E-03 & 1.25E-02 & 2.26E-02 & 6.28E-02 & 0.74  \\  \cline{2-7}
        ~ & aniso: g=0.9 & 8.18E-03 & 1.26E-02 & 2.28E-02 & 6.36E-02 & 0.74  \\ \hline
        \multirow{2}*{bias ($\mathcal{B})$} & iso: g=0 & 1.12E-04 & 2.62E-04 & 6.44E-04 & 6.22E-03 & 1.44  \\  \cline{2-7}
        ~ & aniso: g=0.9  & 1.09E-04 & 2.01E-04 & 3.95E-04 & 5.42E-03 & 1.41 \\ \hline
    \end{tabular}
\end{table}
\subsubsection{Comparison of complexity}
\begin{figure}[ht]
	\centering
\subfigure[isotropic]
{
		\includegraphics[width=0.45\linewidth]{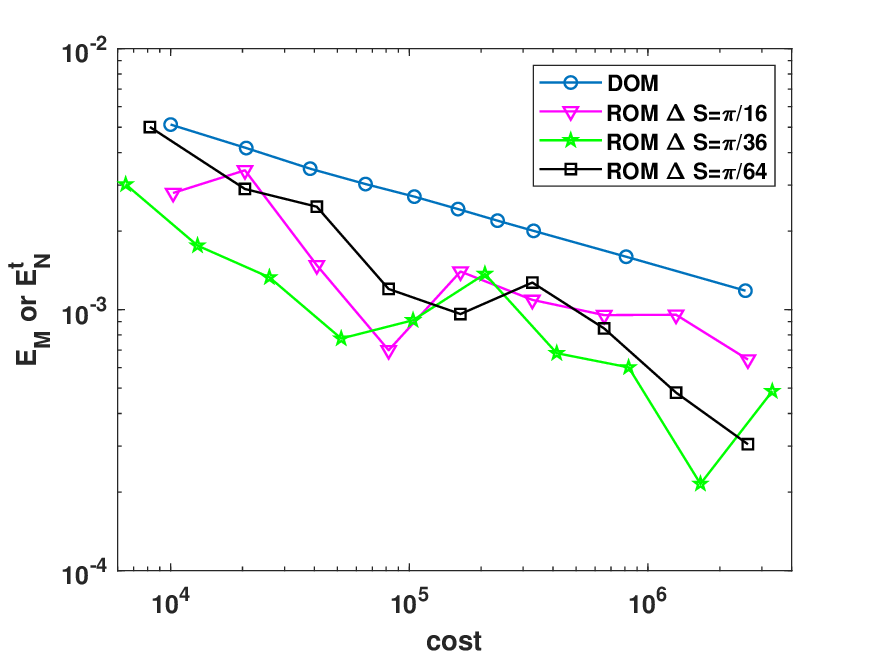}
}
\subfigure[anisotropic]
{
		\includegraphics[width=0.45\linewidth]{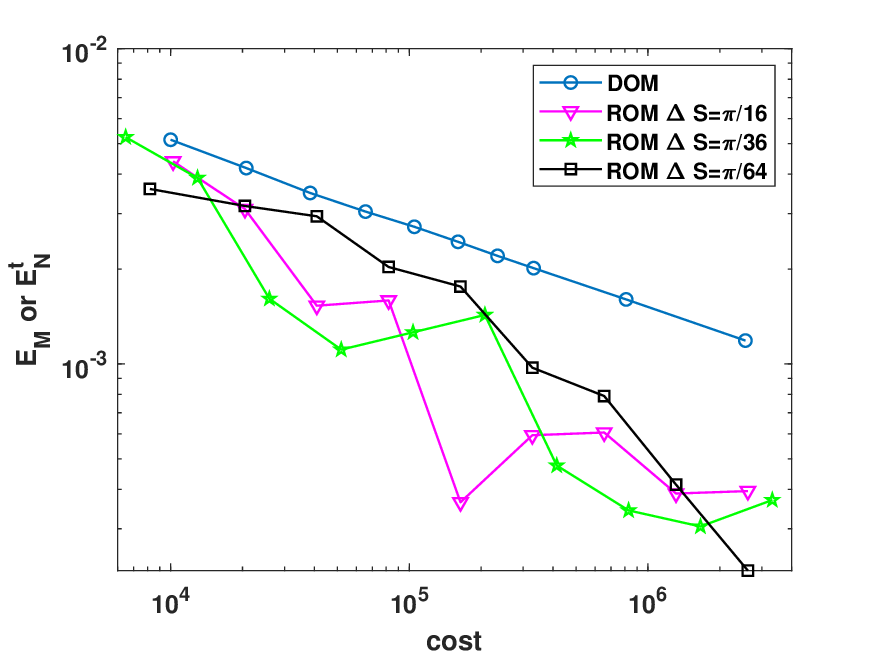}
}

		\caption{Comparison of the DOM and ROM accuracy in terms of the computational cost in X-Y geometry:
For different numbers of ordinates $M = 100, 144, 196, 256, 324, 400, 484, 576, 900, 1600, 1936$, the costs of DOM are assumed to be $M^2$, and the obtained numerical solutions are compared with the reference solution, i.e., the values of $E_M = \|\phi_M - \phi^{ref}\|_2$ are displayed. For ROM, the number of ordinates for each sampled simulation is fixed at $N$, and as the number of sampled simulations $t$ increases from 2 to 20480, the cost is assumed to be $tN^2$. For $N =16,36,64$, the obtained average solution is denoted by $\phi^t_N$, and $E^t_N = \|\phi^t_N - \phi^{ref}\|_2$ are displayed.
 (a):isotropic scattering kernel; (b): anisotropic scattering kernel.
   }
		\label{fig:costvsbias}
\end{figure}

Similar to the slab geometry, we numerically compare the costs of DOM and ROM. In DOM, we compare the obtained numerical solutions with the reference solution calculated using a very fine mesh. Different values of $M = 100, 144, 196, 256, 324, 400, 484, 576, 900, 1600, 1936$ are considered, along with the $L^2$ errors $E_M = \|\phi_M - \phi^{ref}\|_2$. The points $(M^2, E_M)$ are shown in Figure \ref{fig:costvsbias}, where $M^2$ represents the computational cost of solving the $M$ ordinates DOM system. The actual computational cost of solving the DOM system is some constant times $M^2$; the specific value of this constant depends on the employed spatial discretization and mesh size. We use $M^2$ here for simplicity, as when the same spatial discretization and meshes are used in both DOM and ROM, this constant is the same for both methods. For ROM, to consider results obtained at different computational costs, we fix $\Delta S$ and increase the number of sampled simulations $t$. Three different values of $N = 16, 36, 64$ are considered, and the computational cost of ROM is then $tN^2$. The obtained average solution is denoted by $\phi^t_N$, and we define $E^t_N = \|\phi^t_N - \phi^{ref}\|_2$. The points $(tN^2, E^t_N)$ are displayed in Figure \ref{fig:costvsbias}. It is observed that ROM can achieve better accuracy with the same computational cost.

\subsubsection{ Ray effect mitigation}
\begin{figure}[ht]
	\centering
\subfigure[DOM:36 ordinates]
{
		\includegraphics[width=0.45\linewidth]{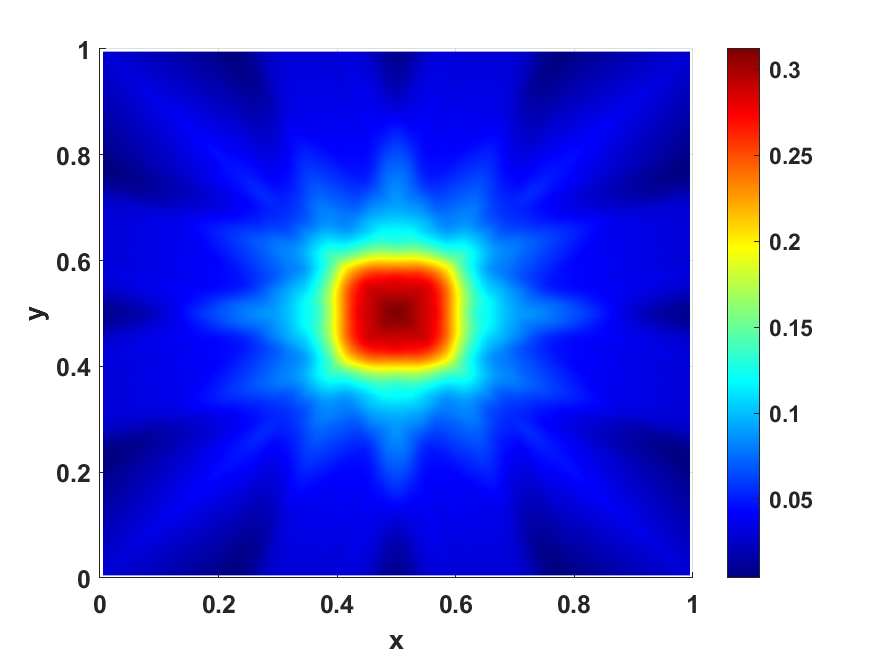}
  \label{fig:L3_ex2}
}
\subfigure[DOM:144 ordinates]
{
		\includegraphics[width=0.45\linewidth]{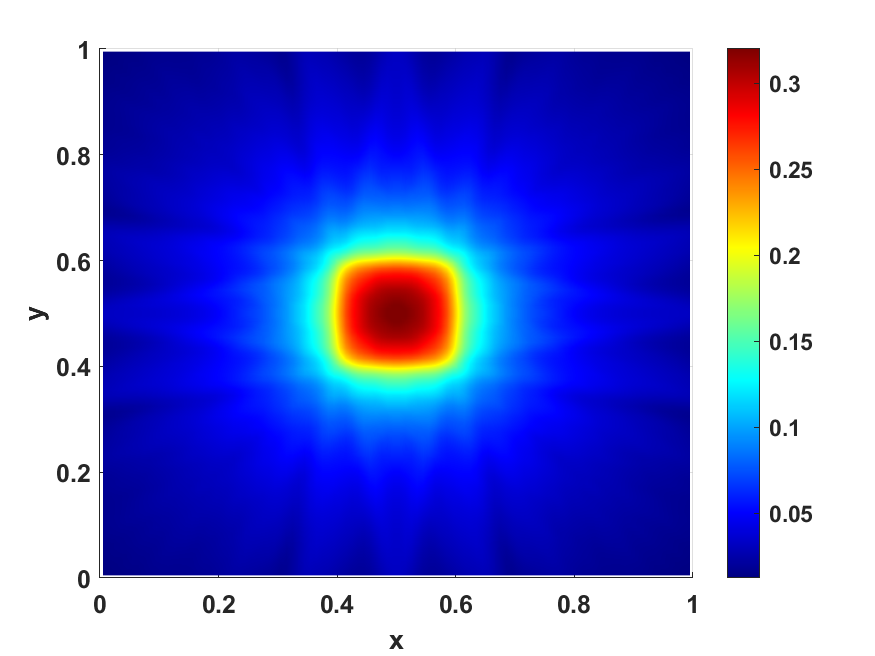}
  \label{fig:L6_ex2}
}

\subfigure[ROM:36 ordinates,16 simulations]
{
		\includegraphics[width=0.45\linewidth]{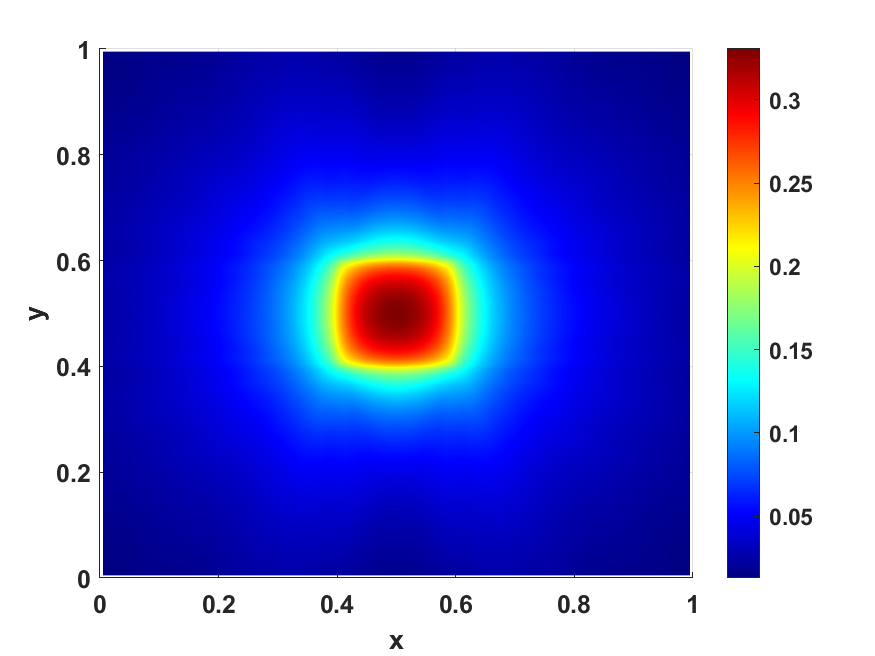}
  \label{fig:ROM_L3_16times_ex2}
}
\subfigure[the solution at $x=0.145$]
{
		\includegraphics[width=0.45\linewidth]{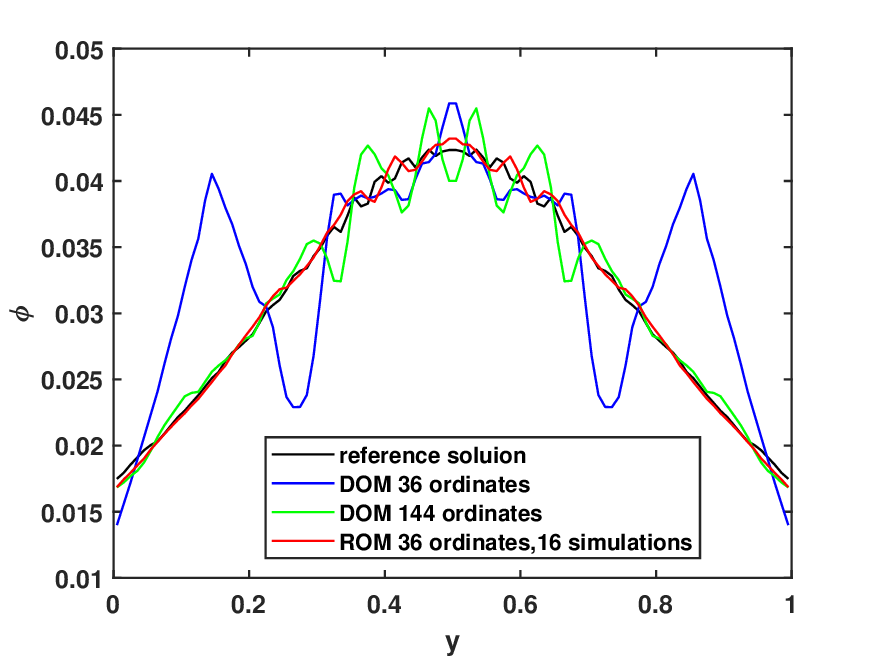}
  \label{fig:ref}
}

		\caption{Example \ref{2D disk_ray effect}. Comparison of the total density obtained by DOM and ROM. (a):the density of DOM with 36 ordinates,(b):the density of DOM with 144 ordinates,(c)the expectation of average density of ROM with 36 ordinatesand 16 simulations, (d): the density profiles using different methods along the $y$ axis at $x=0.15$. The cost of (b) is the same order of magnitude as (c). }
  \label{fig:Comparison of complexity}
\end{figure}

To quantify the ray effect, we define the following metric:
\begin{equation}
    e_{\infty} = \max_{i,j} \left| \phi\left(x_{i+\frac{1}{2}}, y_{j+\frac{1}{2}}\right) - \phi^{ref}\left(x_{i+\frac{1}{2}}, y_{j+\frac{1}{2}}\right) \right|,
\end{equation}
where the reference solution $\phi^{ref}$ is computed by Gaussian quadrature with $80400$ ordinates.
Figure \ref{fig:Comparison of complexity} displays the numerical results for DOM and ROM. Figures \ref{fig:L3_ex2} and \ref{fig:L6_ex2} show the solutions of DOM calculated with $36$ and $144$ ordinates, respectively, and both exhibit a very pronounced ray effect. The results shown in Figure \ref{fig:ROM_L3_16times_ex2} are obtained using ROM with $16$ sampled simulations of $36$ ordinates, which maintains a computational cost comparable to that of DOM with $144$ ordinates as seen in Figure \ref{fig:L6_ex2}. The values of $e_{\infty}$ for the three different tests are as follows: DOM with $36$ ordinates is $0.0373$, DOM with $144$ ordinates is $0.0206$, and ROM with $16$ sampled simulations of $36$ ordinates is $0.012$. This demonstrates that ROM can improve solution accuracy with the same computational cost. As observed in Figure \ref{fig:ref}, the ray effect is also mitigated.

\begin{figure}[ht]
	\centering
\subfigure[4 ordinates,5 simulations]
{
		\includegraphics[width=0.3\linewidth]{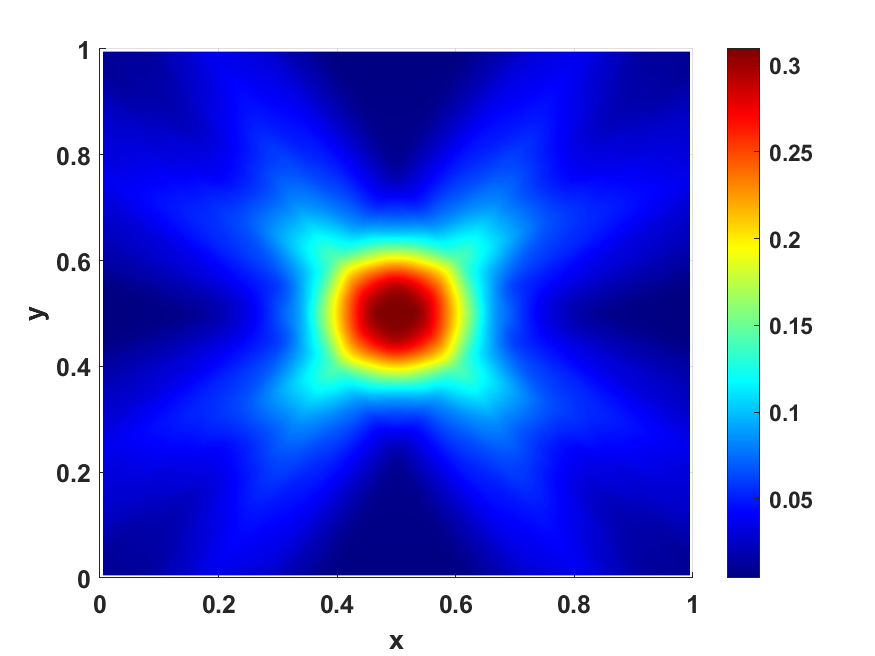}
}
\subfigure[16 ordinates,5 simulations]
{
		\includegraphics[width=0.3\linewidth]{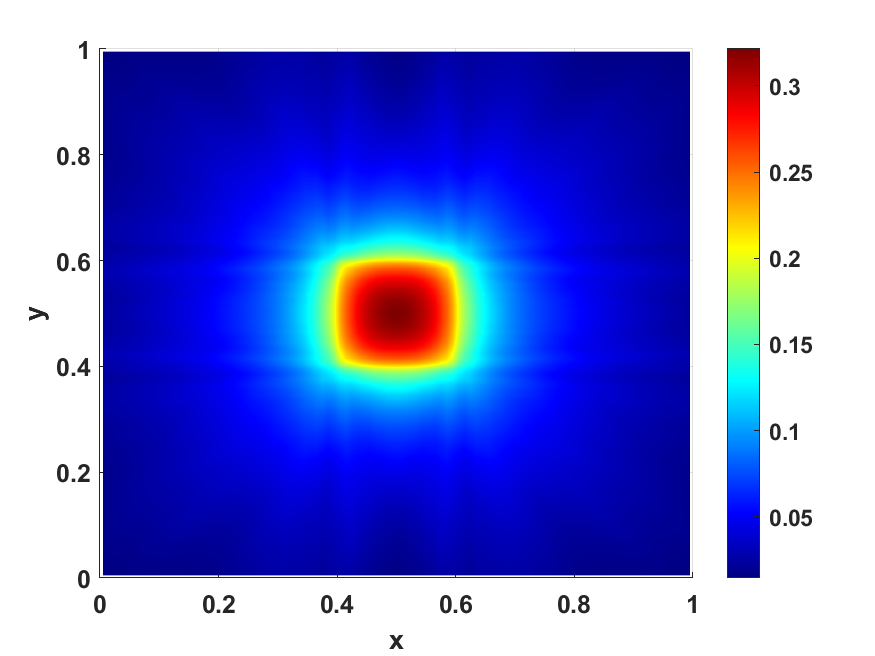}
}
\subfigure[36 ordinates,5 simulations]
{
		\includegraphics[width=0.3\linewidth]{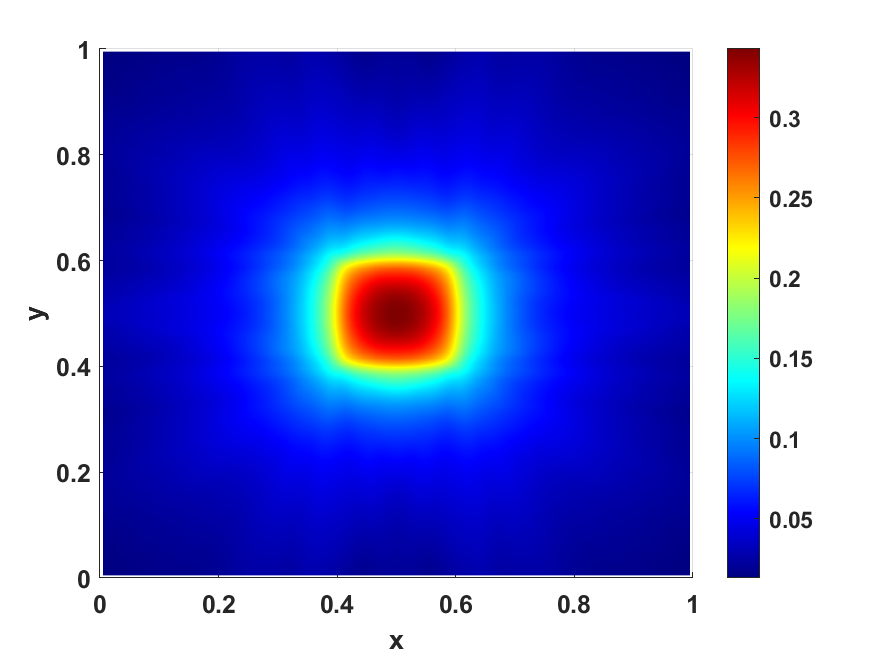}
}

\subfigure[4 ordinates,20 simulations]
{
		\includegraphics[width=0.3\linewidth]{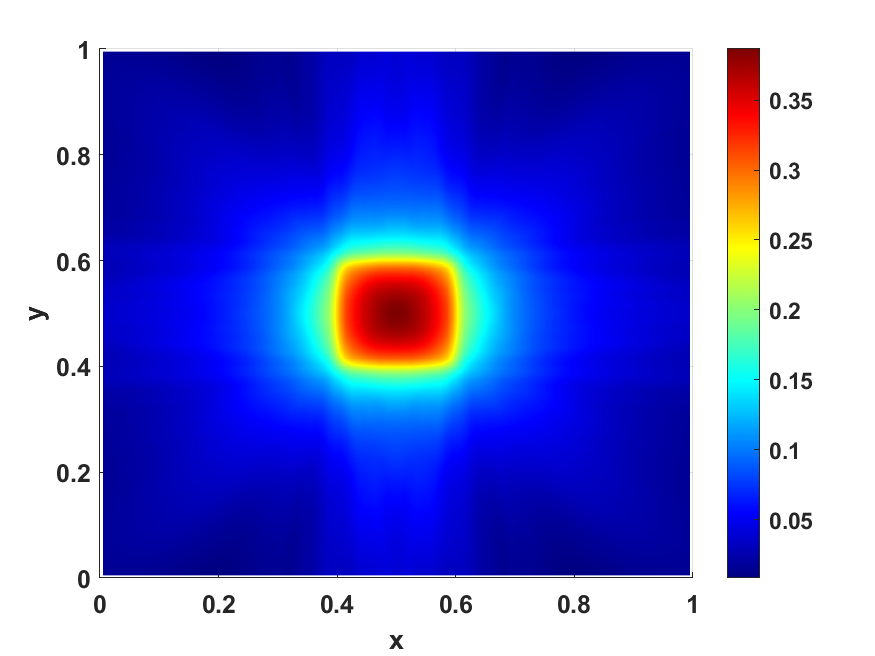}
}
\subfigure[16 ordinates,20 simulations]
{
		\includegraphics[width=0.3\linewidth]{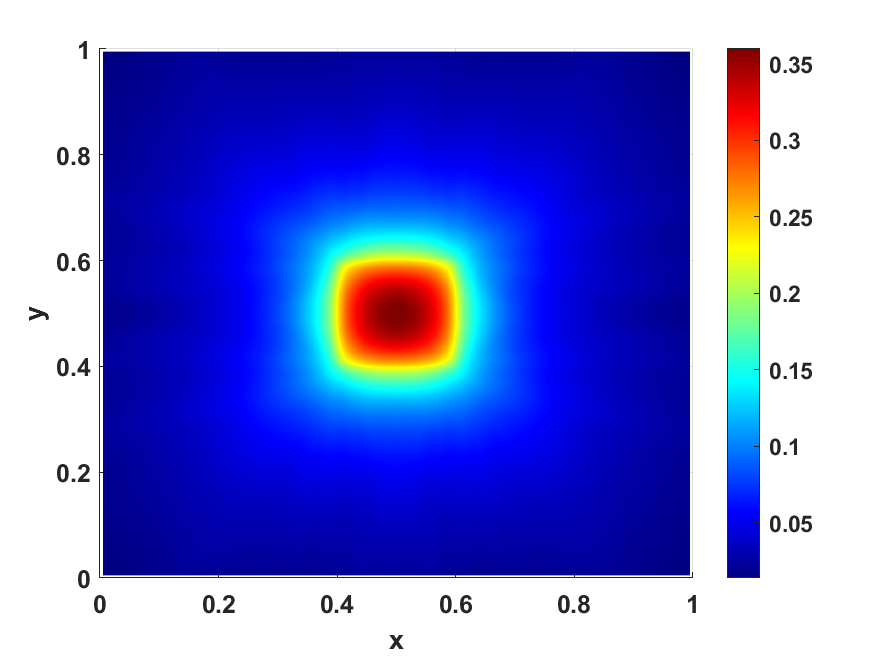}
}
\subfigure[36 ordinates,20 simulations]
{
		\includegraphics[width=0.3\linewidth]{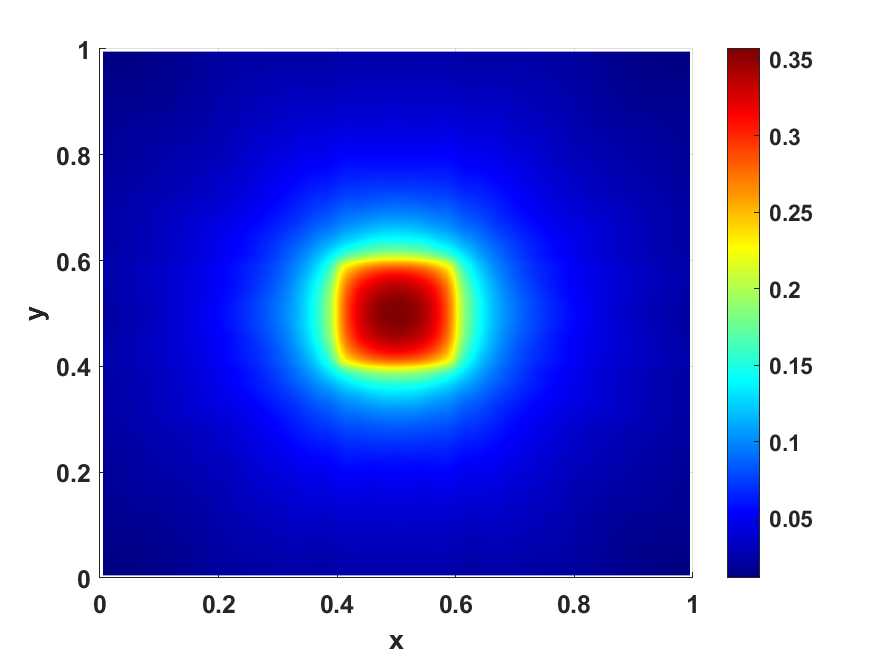}
}

\subfigure[4 ordinates,50 simulations]
{
		\includegraphics[width=0.3\linewidth]{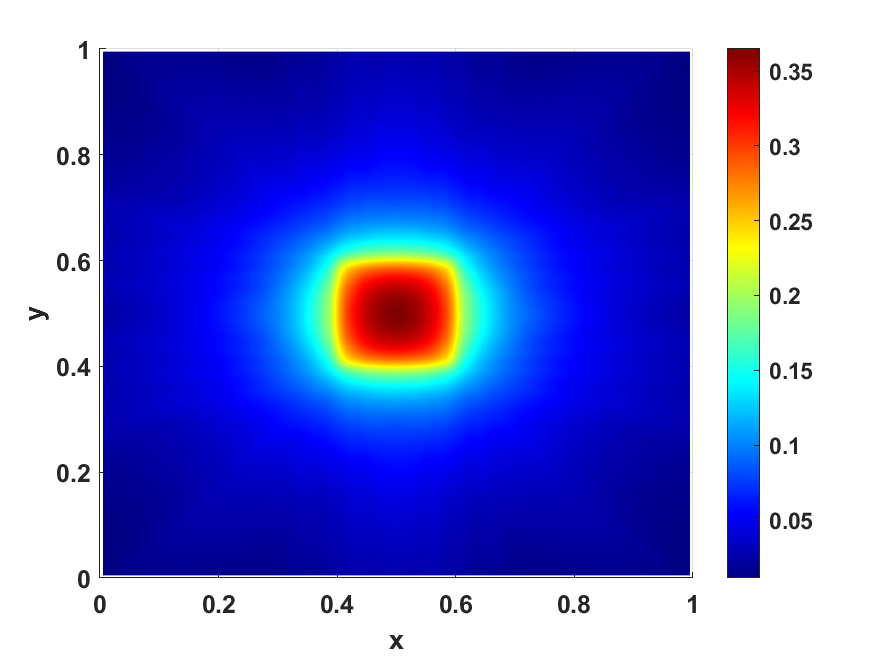}
}
\subfigure[16 ordinates,50 simulations]
{
		\includegraphics[width=0.3\linewidth]{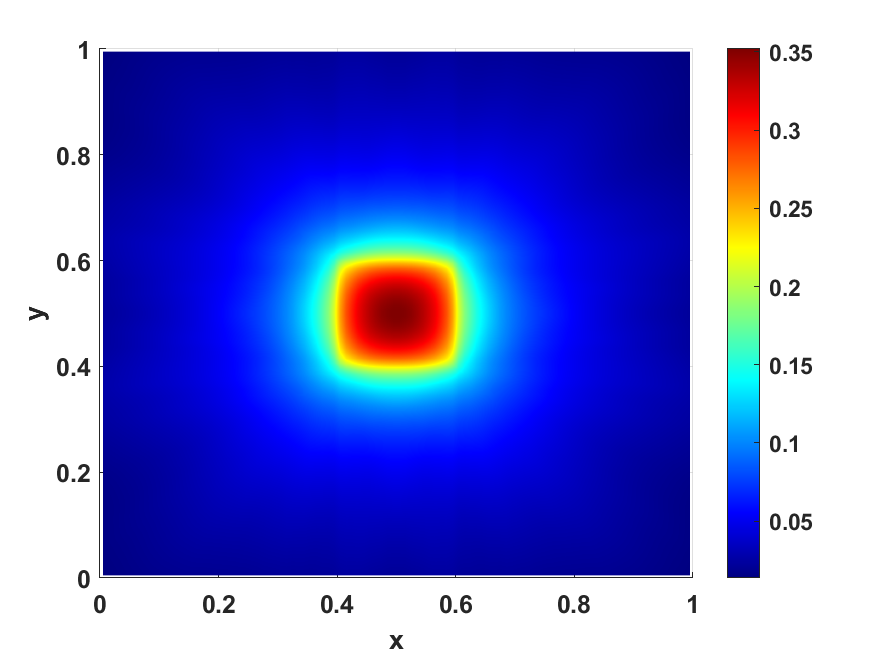}
}
\subfigure[36 ordinates,50 simulations]
{
		\includegraphics[width=0.3\linewidth]{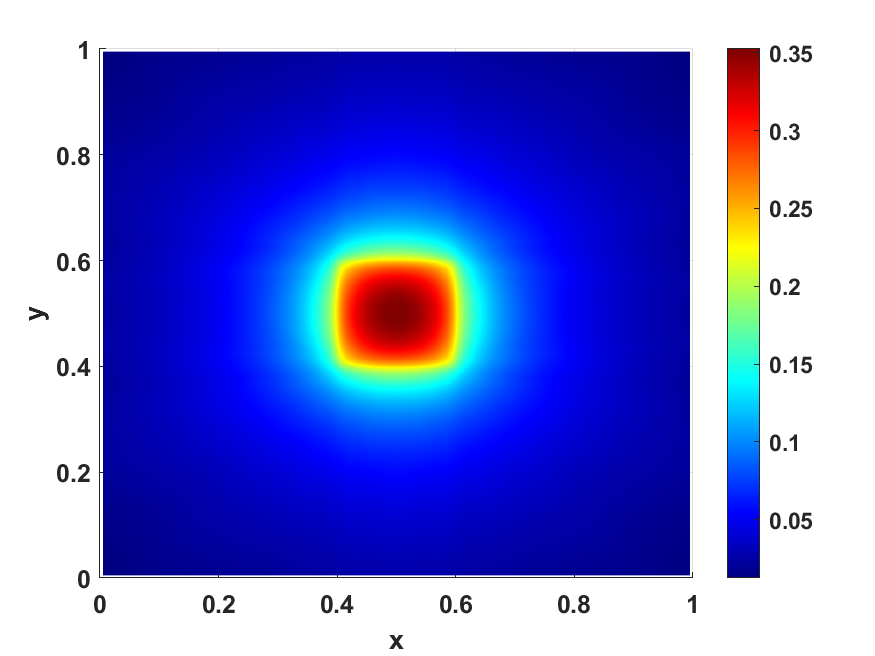}
}
\subfigure[4 ordinates,$x=0.345$]
{
		\includegraphics[width=0.3\linewidth]{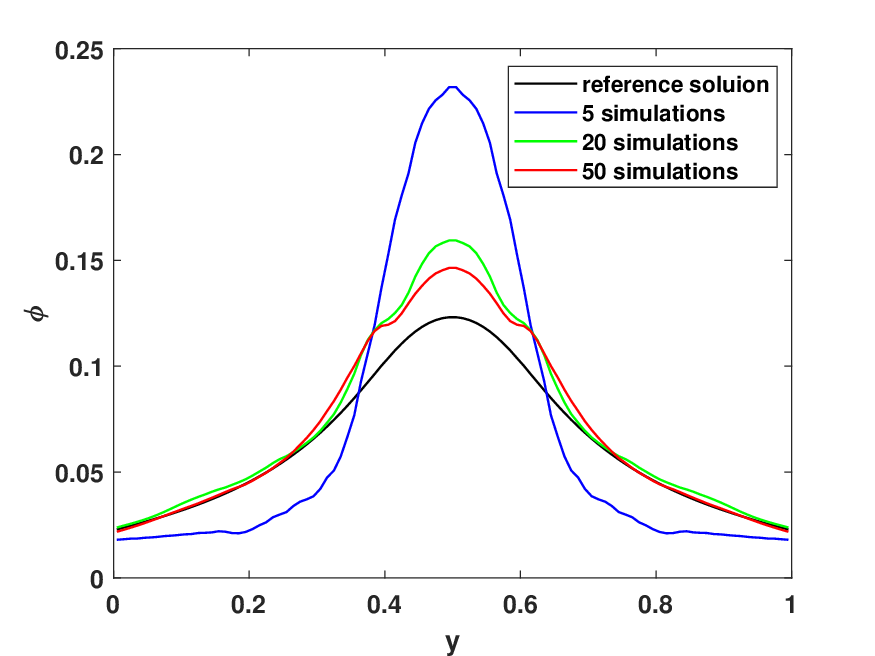}
}
\subfigure[16 ordinates,$x=0.345$]
{
		\includegraphics[width=0.3\linewidth]{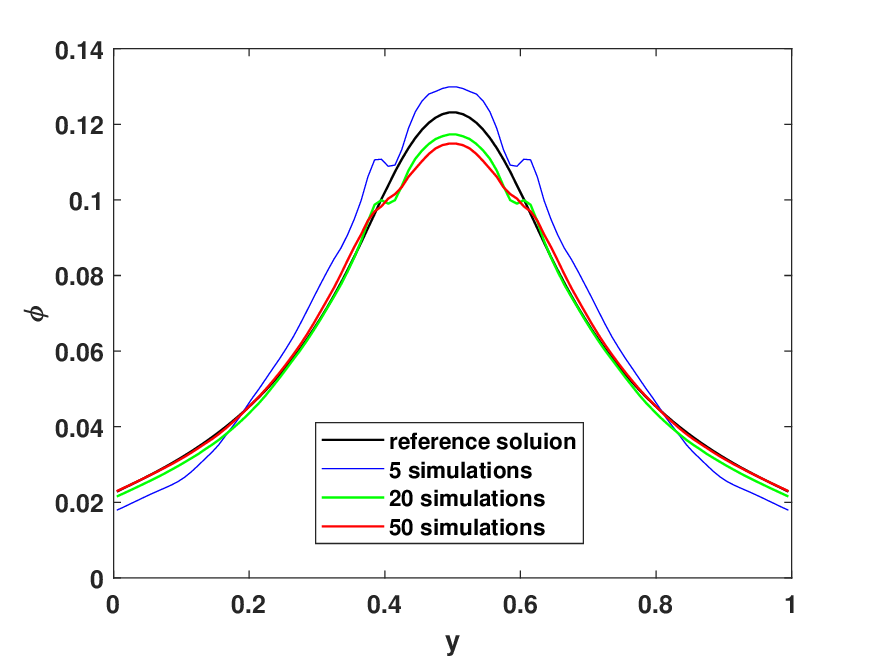}
}
\subfigure[36 ordinates,$x=0.345$]
{
		\includegraphics[width=0.3\linewidth]{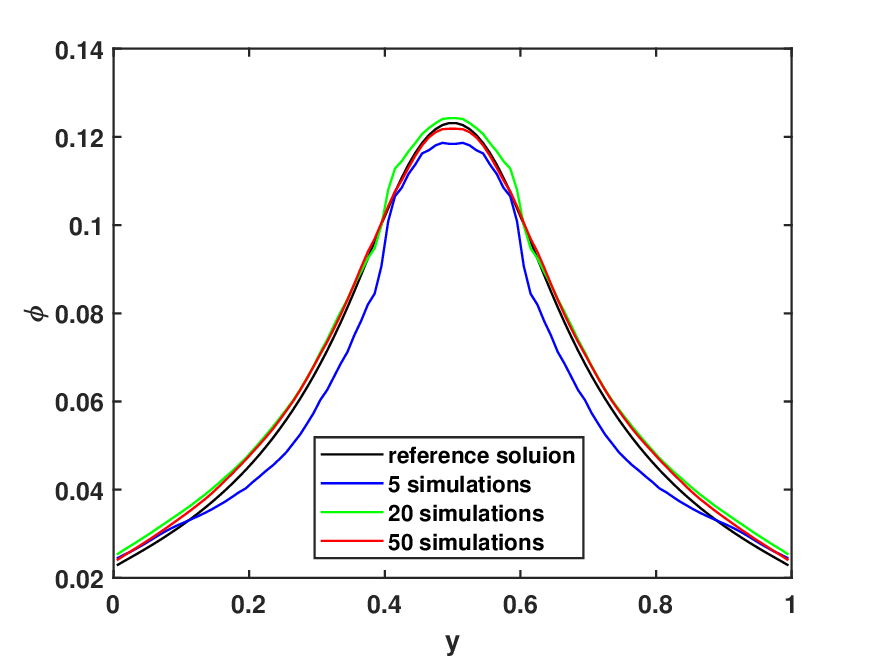}
}

		\caption{Example \ref{2D disk_ray effect}(isotropic). The expectation of the total density of ROM with different number of sampled simulations and ordinates. The fourth row displays the density profiles along the cross-section at $x=0.345$ from the three subplots in the same column.}
		\label{fig:ROMnumerical_5.1_iso}
\end{figure}

Figure \ref{fig:ROMnumerical_5.1_iso} further demonstrates the ability of ROM to mitigate the ray effect. The expectation of the average density calculated by ROM with different number of ordinates and different number of sampled simulations are displayed. We can observe that the ray effect is effectively mitigated after multiple parallel calculations and taking the expectations, even when the sampled simulations use a very small number of ordinates. The fourth row shows the cross-sections at $x=0.345$ from the three subplots above. The solution accuracy improves with an increase in the number of samples. The results with the anisotropic scattering are similar.

\subsubsection{The lattice problem}
\begin{figure}
    \centering
    \includegraphics[width=0.4\linewidth]{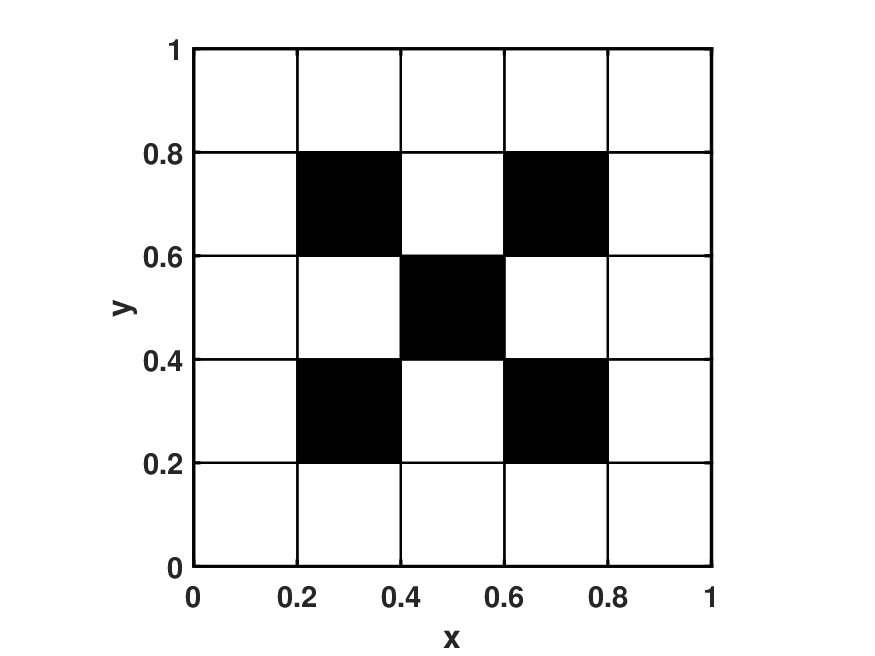}
    \caption{Different source term for lattice problem. The black area is 1, and the other areas are vacant.}
    \label{fig:layout_of_lattice_problem}
\end{figure}

\begin{figure}[ht]
	\centering
\subfigure[4 ordinates]
{
		\includegraphics[width=0.3\linewidth]{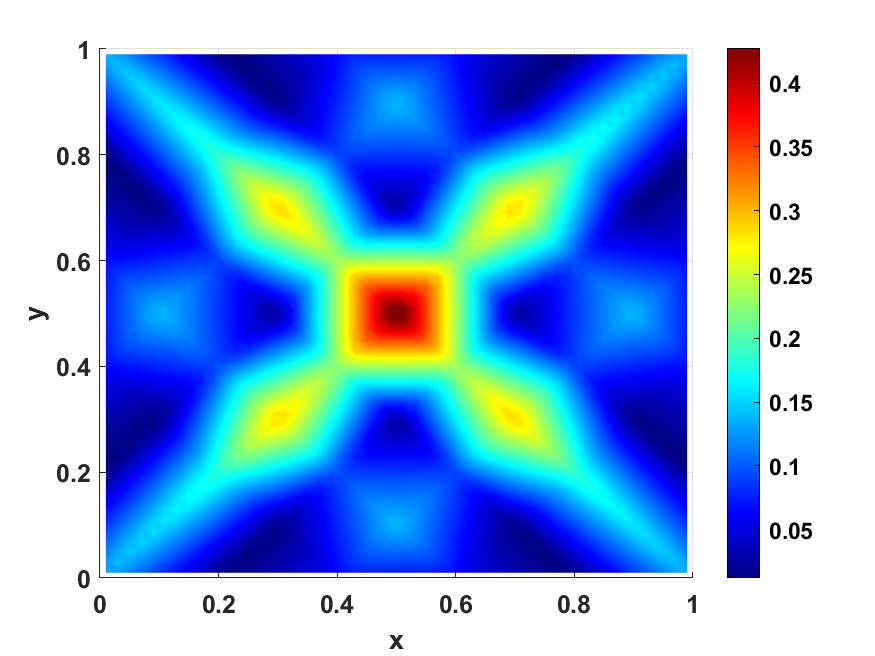}
}
\subfigure[16 ordinates]
{
		\includegraphics[width=0.3\linewidth]{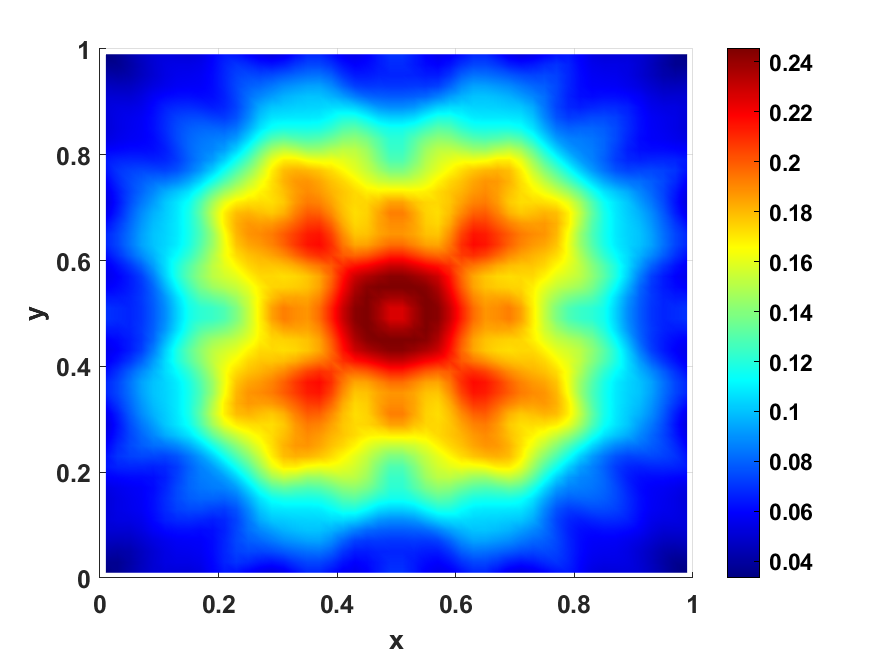}
}
\subfigure[36 ordinates]
{
		\includegraphics[width=0.3\linewidth]{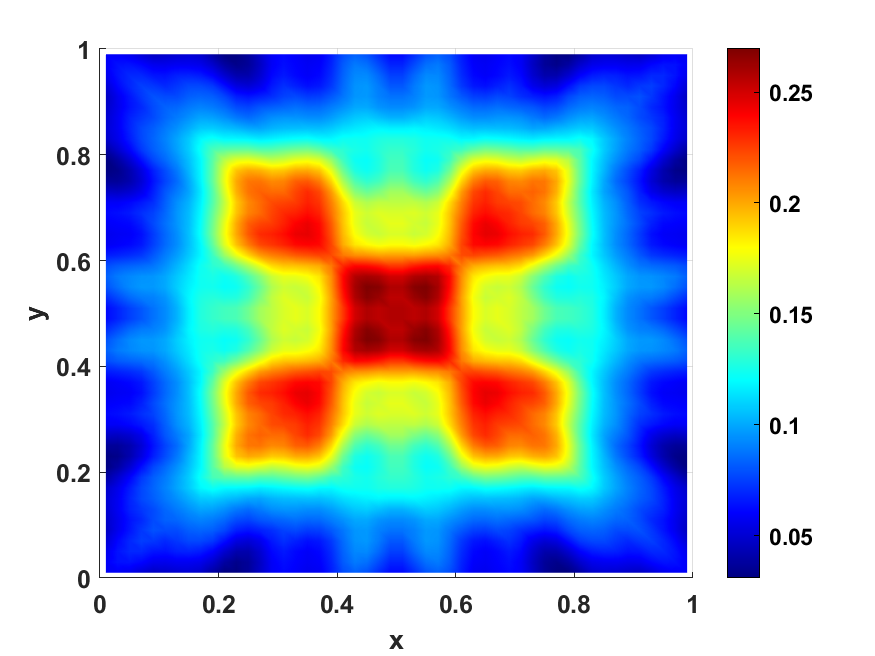}
}

\subfigure[64 ordinates]
{
		\includegraphics[width=0.3\linewidth]{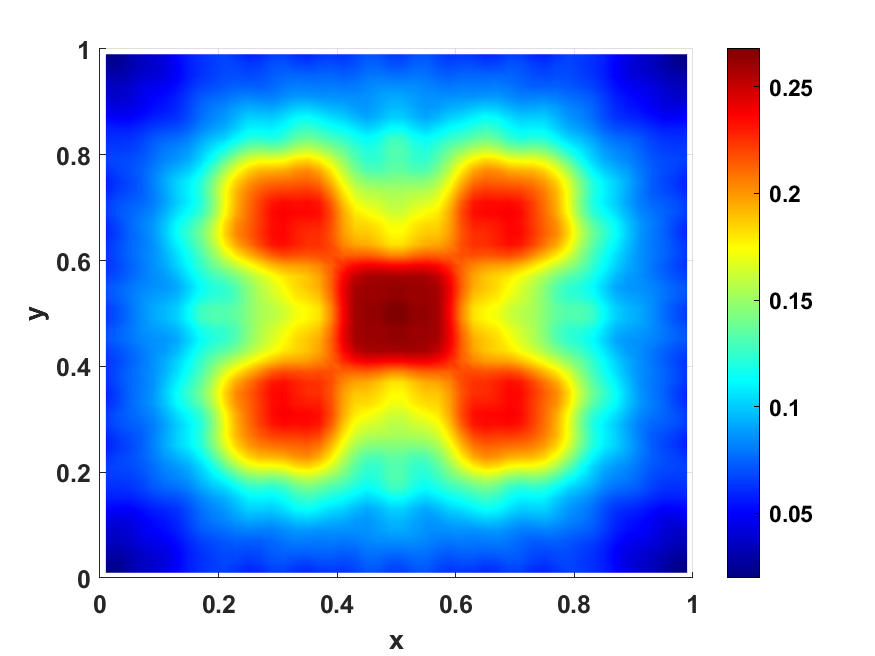}
}
\subfigure[100 ordinates]
{
		\includegraphics[width=0.3\linewidth]{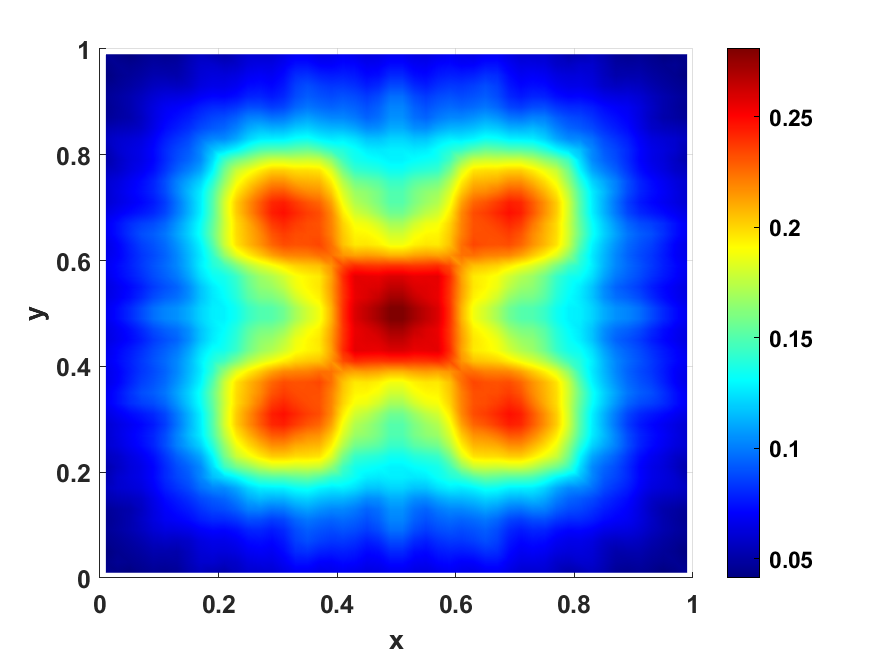}
}
\subfigure[144 ordinates]
{
		\includegraphics[width=0.3\linewidth]{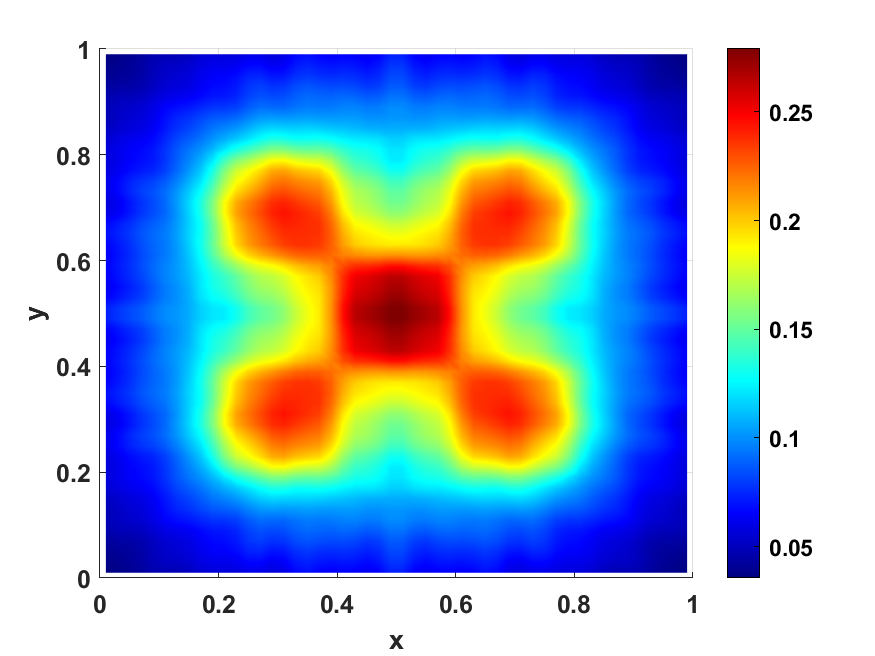}
}

		\caption{The total density obtained by DOM with different number of ordinates of the uniform quadrature. }
		\label{fig:numerical_eliminate_ray_effect}
\end{figure}

This example is to demonstrate that
the ROM can mitigate the ray effect independent of the problem setup. A benchmark test is the lattice problem in X-Y plane. The spatial domain is $(x,y) \in [0,1] \times [0,1]$ and the cross sections are $\sigma_{T}=1$, $\sigma_{S}=0.5$. The layout of the source term $q(x,y)$ is shown in Figure \ref{fig:layout_of_lattice_problem}. The number of spatial cells is $50\times 50$ and Diamond difference method is employed for the spatial discretization. In Figure \ref{fig:numerical_eliminate_ray_effect}, the average densities calculated with $4$, $16$, $36$, $64$, $100$ and $144$ ordinates DOM are displayed. The ray effects can be visually seen even with 144 ordinate. 

\begin{figure}[htp]
	\centering
\subfigure[5 simulations, 4 ordinates]
{
		\includegraphics[width=0.3\linewidth]{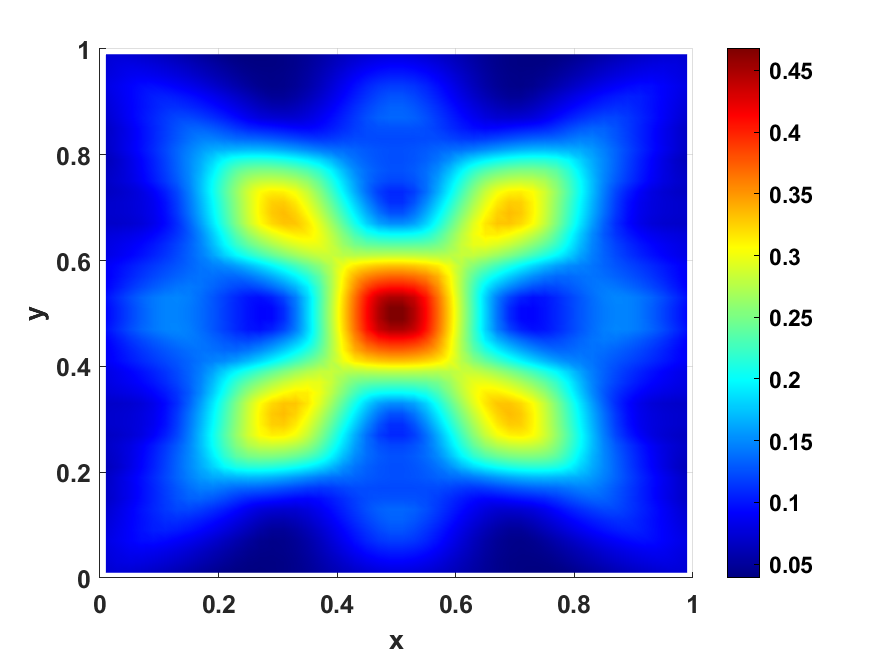}
}
\subfigure[5 simulations, 16 ordinates]
{
		\includegraphics[width=0.3\linewidth]{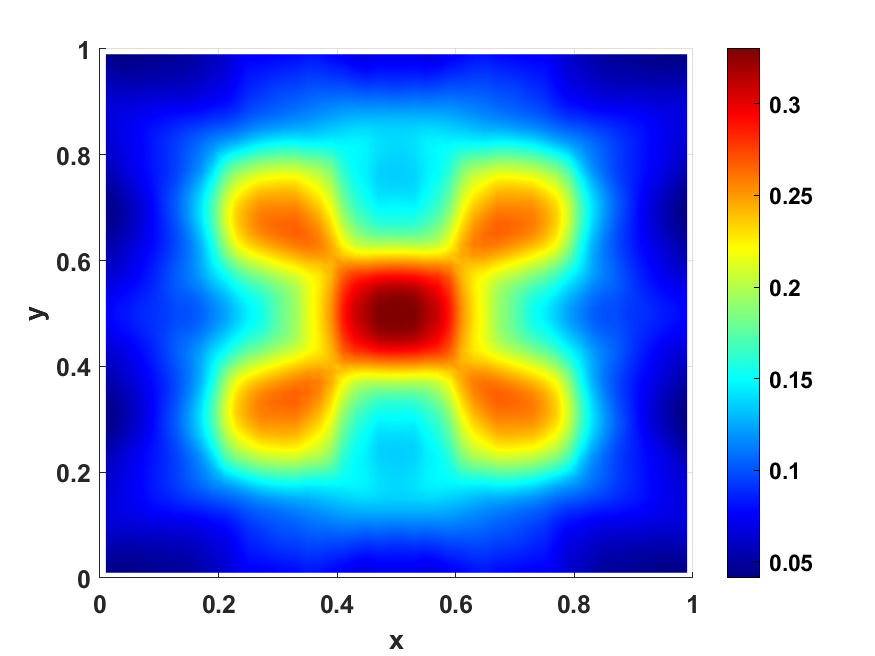}
}
\subfigure[5 simulations, 36 ordinates]
{
		\includegraphics[width=0.3\linewidth]{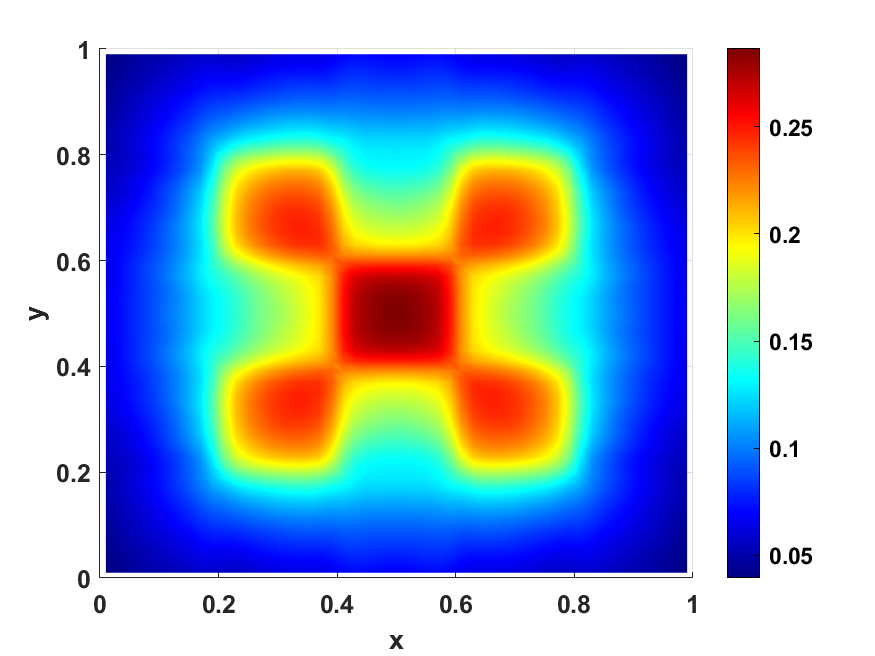}
}

\subfigure[20 simulations, 4 ordinates]
{
		\includegraphics[width=0.3\linewidth]{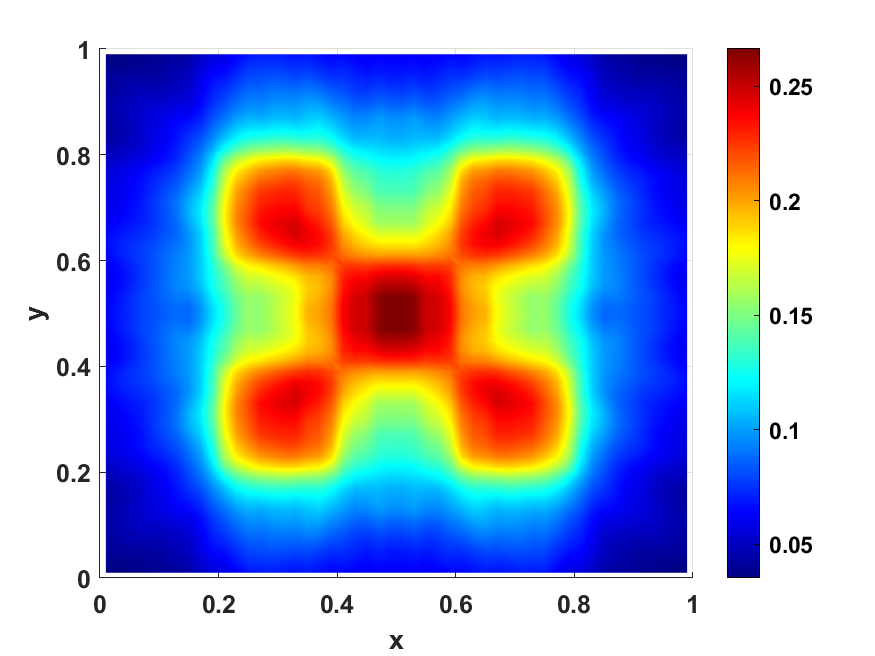}
}
\subfigure[20 simulations, 16 ordinates]
{
		\includegraphics[width=0.3\linewidth]{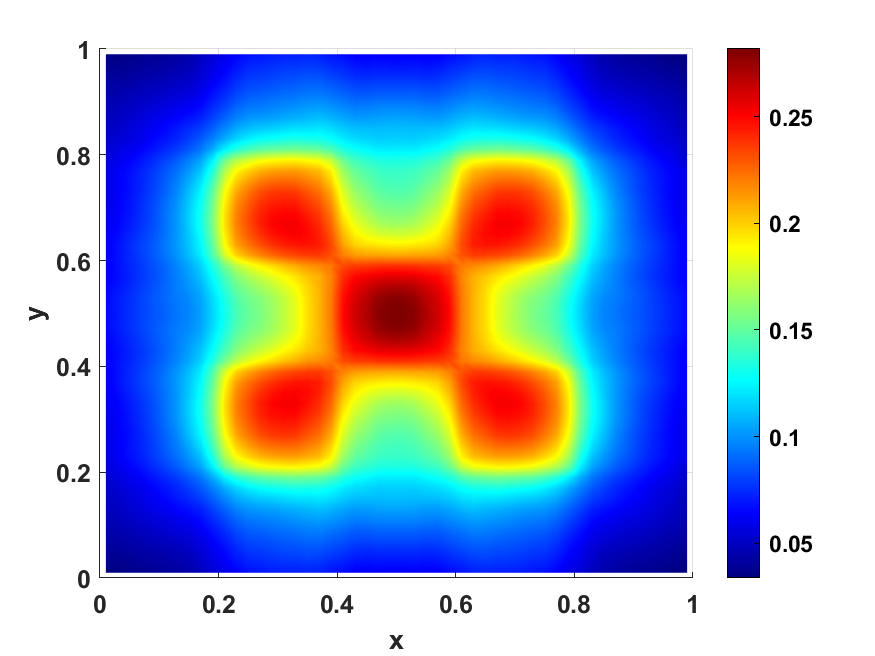}
}
\subfigure[20 simulations, 36 ordinates]
{
		\includegraphics[width=0.3\linewidth]{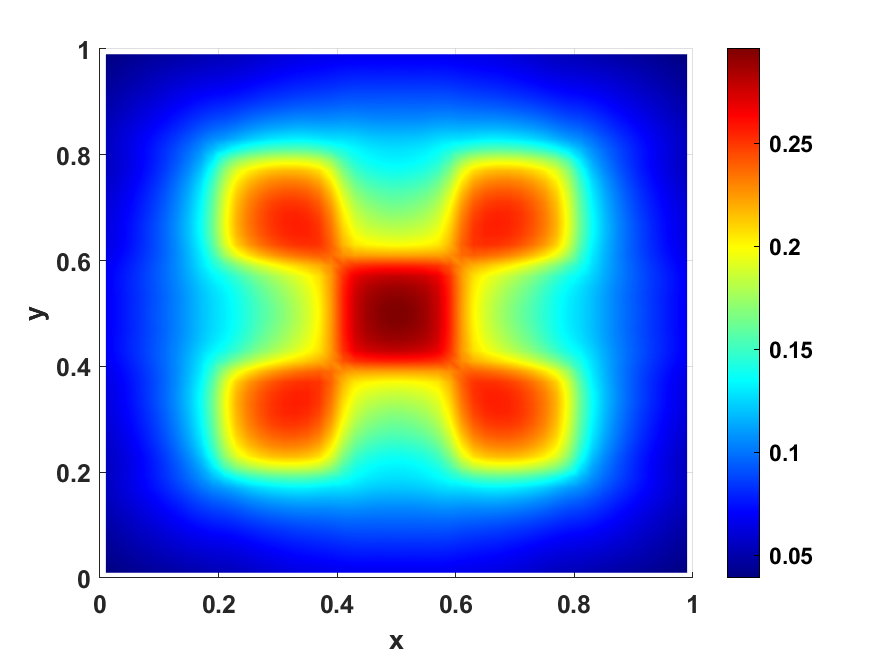}
}
\subfigure[50 simulations, 4 ordinates]
{
		\includegraphics[width=0.3\linewidth]{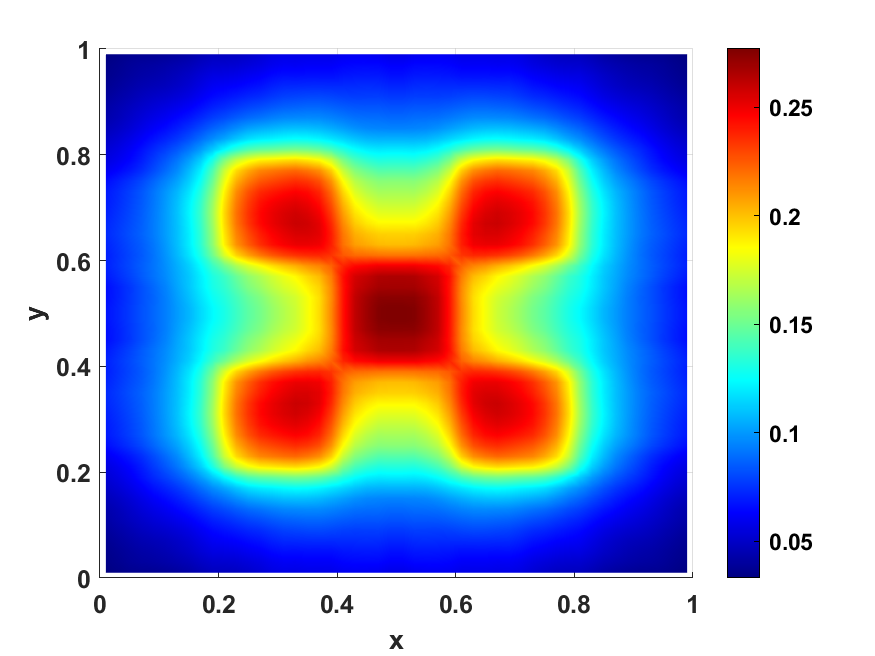}
}
\subfigure[50 simulations, 16 ordinates]
{
		\includegraphics[width=0.3\linewidth]{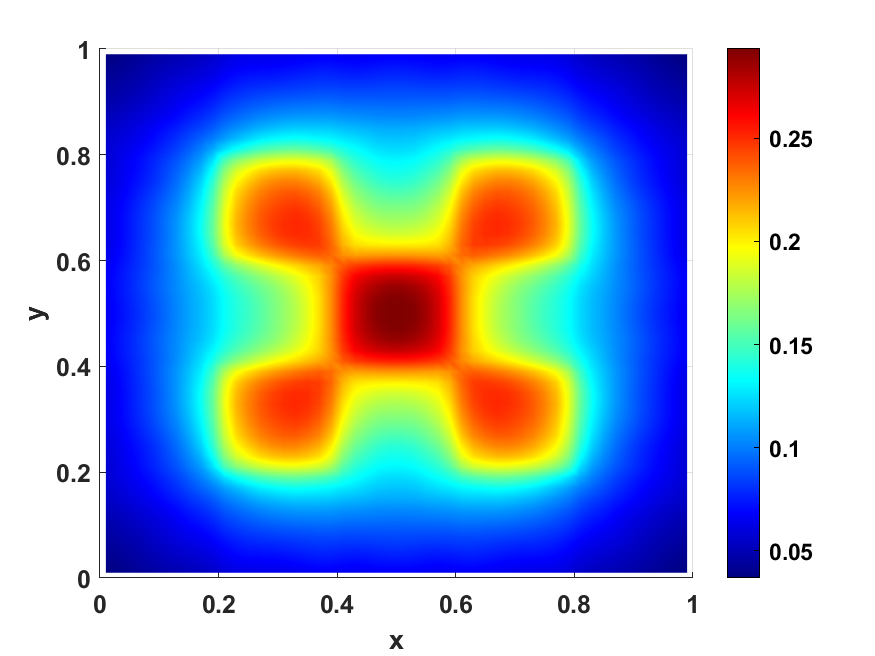}
}
\subfigure[50 simulations, 36 ordinates]
{
		\includegraphics[width=0.3\linewidth]{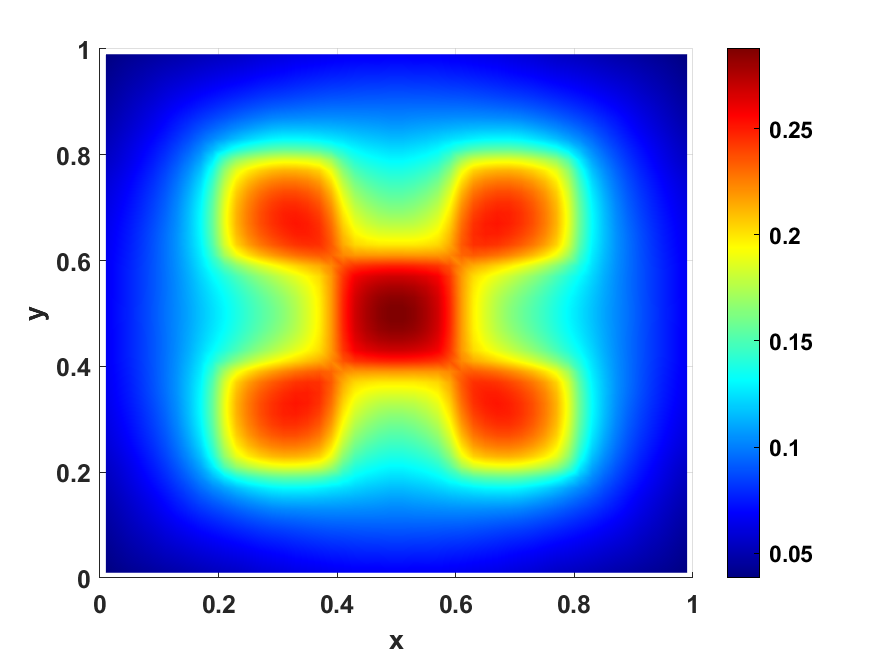}
}

\subfigure[4 simulations, $y=0.29$]
{
		\includegraphics[width=0.3\linewidth]{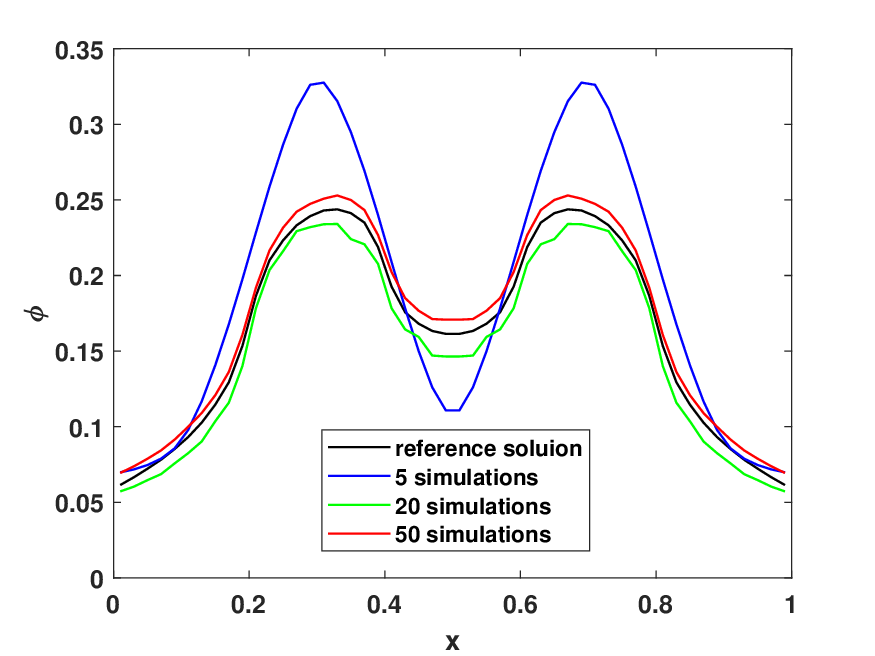}
}
\subfigure[16 ordinates, $y=0.29$]
{
		\includegraphics[width=0.3\linewidth]{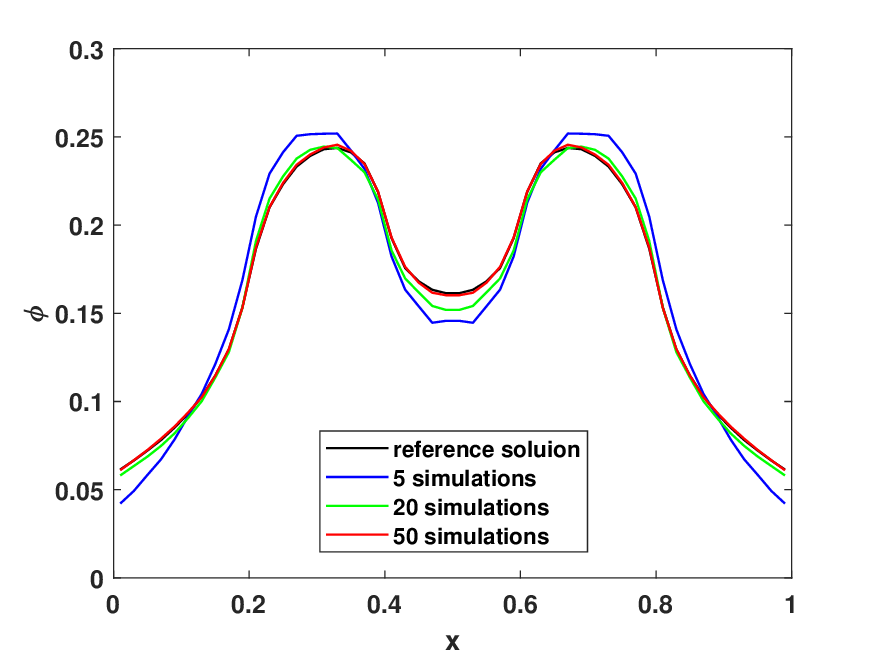}
}
\subfigure[36 ordinates, $y=0.29$]
{
		\includegraphics[width=0.3\linewidth]{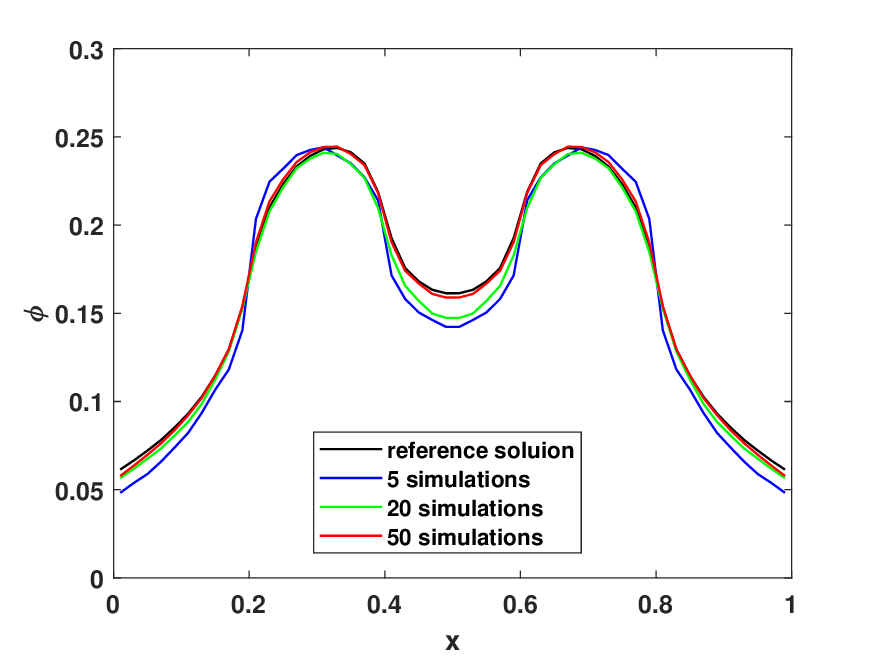}
}

		\caption{The lattice problem. The expectation of the total density obtained by ROM. 
  The fourth row displays the density profiles along the cross-section at $y=0.29$ from the three subplots in the same column.}
		\label{fig:ROMnumerical_eliminate_ray_effect}
\end{figure}

The results of ROM are shown in Figure \ref{fig:ROMnumerical_eliminate_ray_effect}, where the first to third rows show the expectation of average density $\mathbb{E}\phi^{\xi}$ with $5$, $20$ and $50$ simulations, respectively. As can be seen from the results, the ray effect is invisible in the ROM results when $50$ sampled simulations with $4$ ordinates, 20 sampled simulations with $16$ ordinates, and 5 sampled simulations with 36 ordinates. 
The forth row displays the cross-sections at $y=0.29$ from the three subplots above.  The solution accuracy improves with an increase in the number of samples.

\begin{figure}
	\centering
\subfigure[x=0.19]
{
		\includegraphics[width=0.45\linewidth]{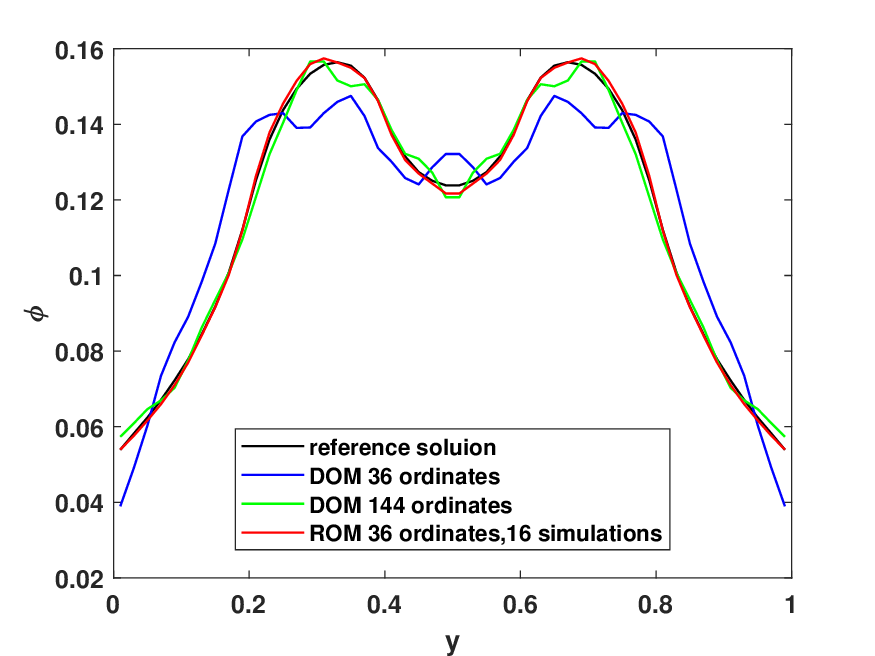}
}
\subfigure[y=0.29]
{
		\includegraphics[width=0.45\linewidth]{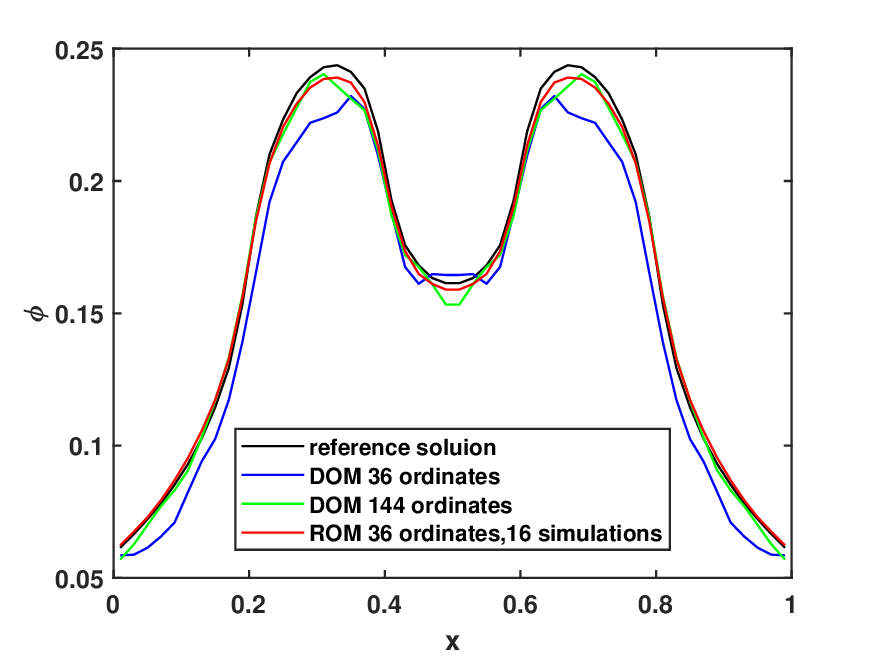}
}

		\caption{ Comparison of the total density obtained by DOM and ROM. (a) the density profiles using different methods along the cross-section at $x=0.19$; (b) the density profiles using different methods along the cross-section at $y=0.29$.
  }
		\label{fig:densityandplotslide}
\end{figure}

Similar to Example \ref{2D disk_ray effect}, we present a comparison of the results obtained by the DOM and ROM with the same computational cost. Figure \ref{fig:densityandplotslide} displays two cuts along both the $x$ and $y$ directions, one at $x=0.19$, the other at $y=0.29$. The four curves represent the reference solution, the solution obtained by DOM with $36$ ordinates, DOM with $144$ ordinates, and ROM with $16$ sampled simulations of $36$ ordinates; the last two have the same computational cost. It can be observed that ROM with $16$ sampled simulations of $36$ ordinates provides the best result and mitigates the ray effect.

\section{Discussion}\label{sec:discussion}
The spatial discontinuities or sharp transitions in the inflow boundary conditions and source terms can induce low regularity in the velocity variable. In slab geometry, discontinuities in velocity are typically at fixed velocity directions, on the other hand, discontinuities in the velocity variable are spatially dependent in X-Y geometry. Consequently, in X-Y geometry, commonly used quadrature sets can only achieve a convergence order around $0.74$, which is the main reason for the observed ray effect.

We have conducted numerical experiments and observed that  the computational cost of the ROM is lower than that of the DOM to achieve the same level of accuracy. The superiority of ROM in alleviating the ray effect is demonstrated through some benchmark tests. Although ROM necessitates multiple runs, which could imply a higher total memory demand, the independence of each run within ROM enables the utilization of the same processor and memory resources across different samples. This feature can lead to a reduced overall memory requirement in comparison to DOM, but there is a price to pay, different samples can no longer run in parallel thus the efficiency reduces. Therefore, one can easily balance memory usage and computational efficiency. 

We remark that other strategies for the random ordinates are possible. For example, $\mathbf{u}_\ell$ does not have to be chosen with uniform probability from $S_{\ell}$. One may make use of the transition kernel for importance sampling of $\mathbf{u}_{\ell}$. If these strategies are adopted, the weights $\omega_{\ell}$ have to be adjusted correspondingly. One may also consider multiscale total and scattering cross sections, in which case, the conservation of mass may become important. In our future work, we will investigate the extension of ROM to multiscale cases, such as the diffusion limit or Fokker-Planck limit cases.




\appendix
\section{Details of the Discrete Ordinate Method}\label{app:detailsDOM}

The steady state RTE reads
\begin{equation}
\b{u} \cdot \nabla \psi(\b{z}, \b{u})+\sigma_{T}(\b{z}) \psi(\b{z}, \b{u})=\sigma_{S}(\b{z})\frac{1}{|S|}\int_{S} P(\b{u'},\b{u})\psi(\b{z}, \b{u'}) \mathrm{d} \b{u'}+q(\b{z}),
\end{equation}
subject to the following inflow boundary conditions:
\begin{equation}\label{eq:appequ_1_bc}
	\psi(\boldsymbol{z}, \boldsymbol{u})=\psi_{\Gamma}^{-}(\boldsymbol{z}, \boldsymbol{u}), \quad \boldsymbol{z} \in \Gamma^{-}=\partial \Omega, \quad \boldsymbol{u} \cdot \boldsymbol{n}_{\boldsymbol{z}}<0.
\end{equation}
After discretizing the velocity with the DOM, equation\eqref{eq:appequ_1_bc} is denoted as
$$
\b{u}_{\ell} \cdot \nabla \psi_{\ell}(\b{z} )+\sigma_{T}(\b{z}) \psi_{\ell}(\b{z})=\sigma_{S}(\b{z}) \sum_{\ell'\in V} w_{\ell'}P_{\ell,\ell'}\psi_{\ell'}(\b{z})+q_{\ell}(\b{z}),\quad \ell\in V,
$$
subject to the following inflow boundary conditions:
$$
	\psi(\b{z}, \b{u}_{\ell})=\psi_{\ell}^{-}(\b{z}, \b{u}), \quad \b{z} \in \Gamma^{-}=\partial \Omega, \quad \b{u}_{\ell} \cdot \b{n}_{\b{z}}<0.
$$

In slab geometry, where $S=[-1, 1]$ and $\Omega=[x_L,x_R]$, the RTE reads, for $(\mu,x)\in[-1,1]\times[x_L,x_R]$:
$$
\mu \partial_x{} \psi(x,\mu)+\sigma_{T}(x) \psi(x,\mu)=\sigma_{S}(x)\frac{1}{2}\int_{-1}^1P(\mu', \mu)\psi(x,\mu')d\mu'+ q(x), 
$$
subject to the boundary conditions:
\begin{equation}\label{eq:1D_bound}
\psi(x_L, \mu)=\psi_L(\mu),\quad \mu>0;\qquad \psi(x_R, \mu)=\psi_R(\mu), \quad\mu<0.
\end{equation}
The DOM in slab geometry takes $V=\{-M,M+1,\cdots,-1,1,\cdots,M-1,M\}$, where $M$ is an integer. The discrete ordinates $\mu_{\ell}$ ($\ell\in V$) satisfy:
$$
0<\mu_{1}<\cdots<\mu_{M-1}<\mu_{M}<1,\qquad \mu_{-\ell}=-\mu_{\ell}.
$$
Therefore, the DOM in slab geometry becomes:
\begin{equation}\label{eq:DOM_1D}
\Big(\mu_{\ell} \partial_{x}+\sigma_{T}(x)\Big) \psi_{\ell}(x)=\sigma_{S}(x) \sum_{\ell' \in V} \omega_{\ell'}P_{\ell,\ell'} \psi_{\ell}(x) + q_{\ell}(x), 
\end{equation}
subject to the boundary conditions: 
$$
	\psi_{\ell}\left(x_{L}\right)=\psi_{L}\left(\mu_{\ell}\right), \quad \mu_{\ell}>0 ; \quad \psi_{\ell}\left(x_{R}\right)=\psi_{R}\left(\mu_{\ell}\right), \quad \mu_{\ell}<0.
$$

Two commonly used quadratures are the Uniform quadrature and the Gaussian quadrature. For the uniform quadrature, $[-1, 1]$ is divided into $2M$ equally spaced cells, each of size $\Delta \mu = 1/M$ when M  is an integer. The values $\{\mu_{\ell} | \ell \in V\}$ represent the midpoint of each cell, i.e., $\mu_{\ell}=\frac{2\ell-1}{2M}$ for $\ell>0$ and $\mu_{\ell}=\frac{2{\ell}+1}{2M}$ for $\ell<0$, while $\omega_{\ell}=1/M$. For the Gaussian quadrature, $\{\mu_{\ell} | \ell \in V\}$ consists of $2M$ distinct roots of Legendre polynomials of degree $2M$ denoted by $L_{2M}(x)$, and the weights $\omega_{\ell}=2/[(1-\mu_{\ell}^2)(L_{2 M}^{\prime}(\mu_{\ell}))^2]$.

For RTE in the X-Y geometry, where the 3D velocity on a unit sphere is projected to a 2D disk, suppose that the DOM has $M$ points in each quadrant. Then, the 2D discrete velocity directions defined on a disk are $\mathbf{\b{u}_{\ell}}=(c_{\ell},s_{\ell})$ for ${\ell}\in \bar V=\{1,2,\cdots M,\cdots,4M\}$ when $M$ is an integer. The DOM in X-Y geometry is expressed as follows, for ${\ell}\in \bar V$ and $(x,y)\in[x_L,x_R]\times[y_B,y_T]$:
$$
\Big(c_{\ell} \partial_{x}+  s_{\ell} \partial_{y}+\sigma_T(x,y)\Big) \psi_m(x,y)
=\sigma_S(x,y) \sum_ {\ell'\in \bar{V}}\bar{\omega}_{\ell'}P_{\ell,\ell'}\psi_{\ell'}(x, y) +q_{\ell}(x,y),
$$
where
$$
(c_{\ell},s_{\ell})=\Big(\left(1-\zeta_{\ell}^2\right)^{\frac{1}{2}} \cos \theta_{\ell},  \left(1-\zeta_{\ell}^2\right)^{\frac{1}{2}} \sin \theta_{\ell}\Big), \quad\mbox{with } \zeta_{\ell} \in (0,1), \theta_{\ell} \in (0,2 \pi).
$$
The boundary conditions become 
\begin{equation*}
\left\{\begin{array}{llll}
\psi_{\ell}\left(x_L, y\right)=\psi_{L, \ell}(y), & c_{\ell}>0 ; & \psi_{\ell}\left(x_R, y\right)=\psi_{R, \ell}(y), & c_{\ell}<0; \\
\psi_{\ell}\left(x, y_B\right)=\psi_{B, \ell}(x), & s_{\ell}>0 ; & \psi_{\ell}\left(x, y_T\right)=\psi_{T, \ell}(x), & s_{\ell}<0.
\end{array}\right.
\end{equation*}

Two kinds of quadratures discussed in \cite{lewis1984computational} are considered. Each quadrant has $M$ discrete velocities, and we only show the details for the first quadrant such that $\zeta_{\ell}\in(0,1)$ and $\theta_{\ell}\in (0,\frac{\pi}{2})$. The discrete velocities in other quadrants are obtained by symmetry.

The first one is referred to as "2D uniform quadrature". Each quadrant has $M=N^2$ ordinates, and the nodes are uniform in the $(\zeta,\theta)$ plane. More precisely, $(\zeta_i,\theta_j)=(\frac{2N-2i+1}{2N},\frac{2j-1}{4N}\pi)$ for $i=1,\cdots,N$; $j=1,\cdots,N$ and $\bar{\omega}_{\ell}=\frac{1}{4N^2}$. Then for $M=N^2$ and all $\ell\in\{1,\cdots, M\}$, there exists a pair of integers $(i,j)$ with $i,j\in\{1,\cdots,N\}$ such that $\ell=(i-1)N+j$ and 
\begin{equation}
(c_{\ell},s_{\ell},\bar{\omega}_{\ell})=\Big(\big(1-\zeta_i^2\big)^{\frac{1}{2}}\cos\theta_j,\big(1-\zeta_i^2\big)^{\frac{1}{2}}\sin\theta_j,\frac{1}{4N^2}\Big).
\label{eq:uniform2D}
\end{equation}

\begin{figure}[htbp]
	\centering
\subfigure[]
{
		\includegraphics[width=0.4\linewidth]{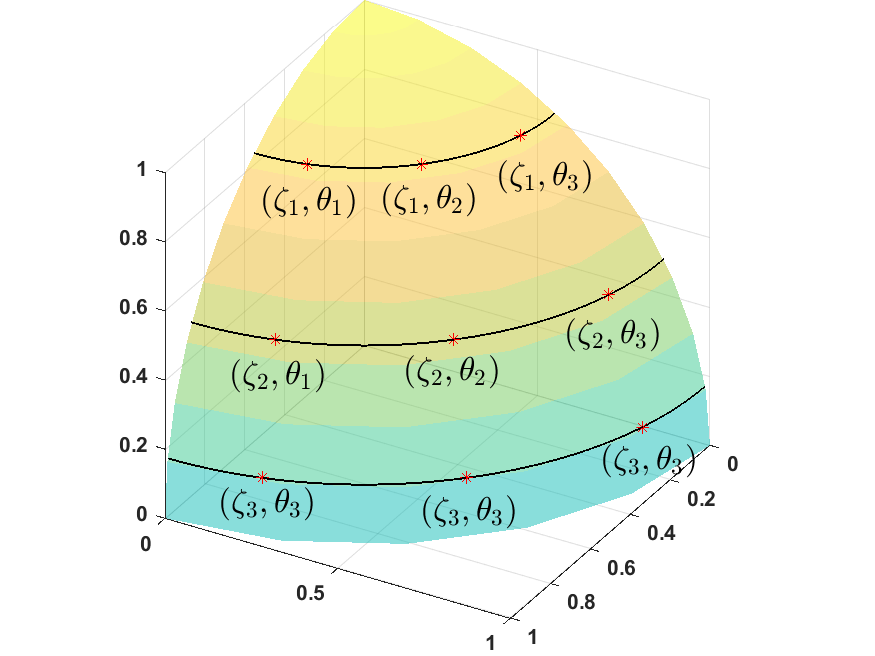}
}
\subfigure[]
{
		\includegraphics[width=0.4\linewidth]{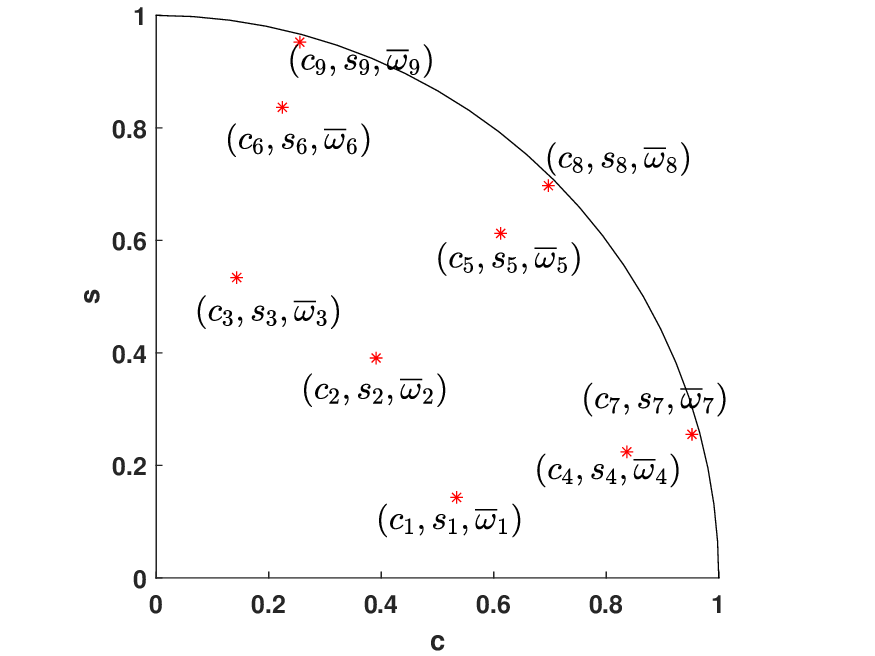}
}
\subfigure[]
{
		\includegraphics[width=0.4\linewidth]{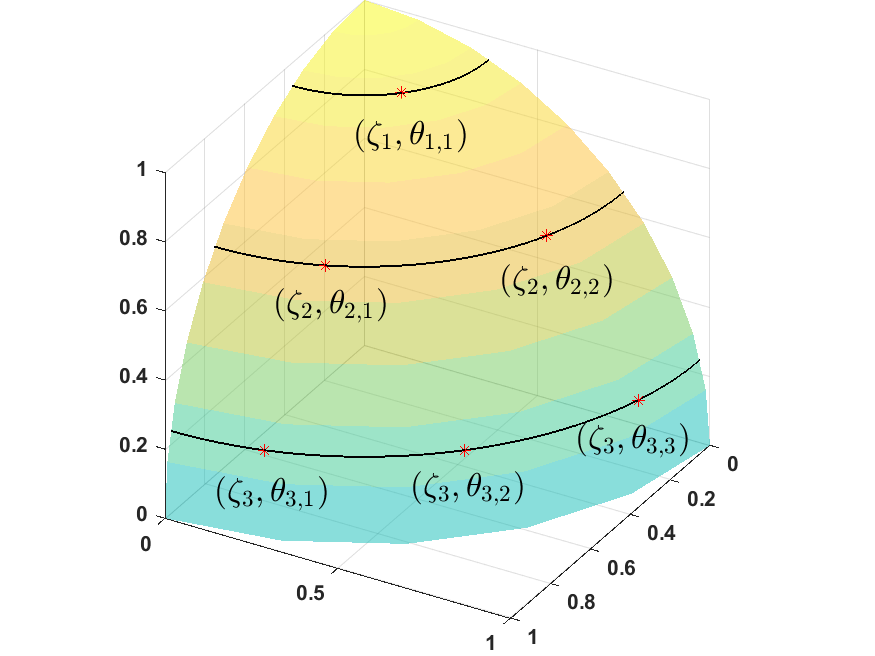}
}
\subfigure[]
{
		\includegraphics[width=0.4\linewidth]{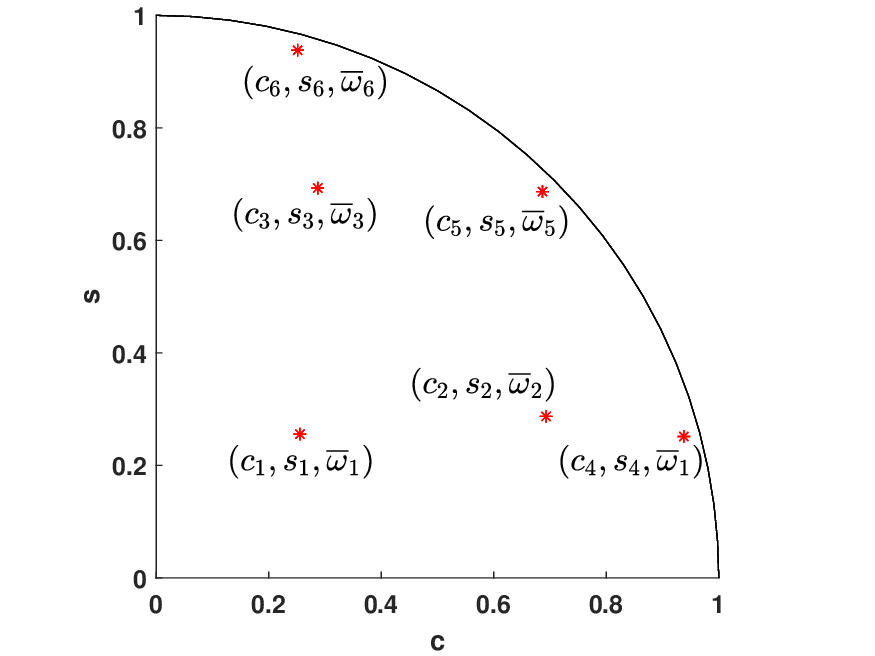}
}
		\caption{Schematic diagram of selected ordinates on the surface of a 3D unit sphere and their corresponding projection to the 2D unit disk. (a)(b) Uniform quadrature; (c)(d) Gaussian quadrature.}
		\label{fig:quadrature_3d}
\end{figure}

The second one is the 2D Gaussian quadrature described in \cite{lewis1984computational}, for which each quadrant has $M=N(N+1) / 2$ ordinates.
Each quadrant has $N$ distinct $\zeta_i$, $i\in\{1,\cdots, N\}$, which are the $N$ positive roots of $L_{2 N}(\zeta)$, the Legendre polynomial of degree $2N$. They are arranged as 

$$
1>\zeta_1>\zeta_2>\cdots>\zeta_N>0.
$$
Each $\zeta_i$ corresponds to $N$ distinct $\theta_{i, j}=\frac{2 j-1}{4 i} \pi, j=1,2, \cdots, i$, and the weight for the velocity direction $(\zeta_i,\theta_{i,j})$ is uniform in $j$ such that
$$
\bar{\omega}_i=\frac{1}{2i\left(1-\zeta_i^2\right)\left[L_{2 N}^{\prime}\left(\zeta_i\right)\right]^2} .
$$
Then for $M=N(N+1) / 2$ and all $\ell\in\{1,\cdots,M\}$, there exists a pair of integers $(i,j)$ such that $i\in\{1,\cdots, N\}$, $1\leq j\leq i$, and $\ell=\frac{i(i-1)}{2}+j$ and

$$
\left(c_\ell,s_\ell, \bar{\omega}_{\ell}\right)=\Big(\big(1-\zeta_i^2\big)^{\frac{1}{2}}\cos\theta_{i,j},\big(1-\zeta_i^2\big)^{\frac{1}{2}}\sin\theta_{i,j},\frac{1}{2i\left(1-\zeta_i^2\right)\left[L_{2N}^{\prime}\left(\zeta_i\right)\right]^2}\Big).
$$
The rest part of the quadrature set can be constructed by symmetry:
$$\bar{\omega}_{\ell}=\bar{\omega}_{\ell+M}=\bar{\omega}_{\ell+2 M}=\bar{\omega}_{\ell+4 M},$$
$$\theta_{\ell}=\pi-\theta_{\ell+M}=\theta_{\ell+2 M}+\pi=-\theta_{\ell+4 M},  \qquad
 \zeta_{\ell}=\zeta_{\ell+M}=\zeta_{\ell+2 M}=\zeta_{\ell+4 M}.
$$
Let's take $N=3$, the chosen discrete ordinates of Uniform and Gaussian quadratures in one quadrant are plotted in Figure \ref{fig:quadrature_3d}. The selected ordinates on the surface of a 3D unit sphere and their corresponding projections to the 2D unit disk are displayed.

\bibliographystyle{siamplain}
\bibliography{references}
\end{document}